\documentclass[letterpaper, 10 pt]{IEEEtran}
 \usepackage{pdfsync}%switch back and forth between the PDF document and the text editor
 \usepackage{braket,amsfonts}
\usepackage{bm,array}
\usepackage{placeins,graphicx,diagbox }
\usepackage{etex,amsmath,amssymb}  % assumes amsmath package installed
\usepackage{pgfplots,stfloats}

\usepackage{amsmath,algorithm} % assumes amsmath package installed
\usepackage{amssymb}  % assumes amsmath package installed
\usepackage{color}
\usepackage{cite}
\newcounter{algo}

\newtheorem{theorem}{Theorem}
\newtheorem{lemma}{Lemma}

\newtheorem{assumption}{Assumption}

\newtheorem{corollary}{Corollary}
\newtheorem{remark}{Remark}

\usepackage{algorithmic}

\usepackage{graphicx,epstopdf}
\usepackage{amsopn}

\usepackage{xspace}
\usepackage{bold-extra}

\usepackage{placeins,diagbox }

\usepackage{wrapfig} % for wrapping the text around the fig
\usepackage{tikz,pgfplots,mathtools,float}
\usepackage{booktabs,caption}
\usepackage[flushleft]{threeparttable}
\usepackage{subfigure}
\usepackage{multirow}
\usepackage{amssymb,latexsym}  % assumes amsmath package installed
\usepackage{extarrows}
\usepackage{dsfont} % for double stroke font
\usepackage{amssymb}  % assumes amsmath package installed
\usepackage{color}
\usepackage{cite}

\colorlet{texcscolor}{blue!50!black}
\colorlet{texemcolor}{red!70!black}
\colorlet{texpreamble}{red!70!black}
\colorlet{codebackground}{black!25!white!25}

\def\Real{\mathbb{R}}

\newcommand{\blue}[1]{{\color{black}#1}}
\newcommand{\red}[1]{{\color{black}#1}}

\title{ Distributed Variable Sample-size Stochastic  Optimization with Fixed Step-sizes  }

\author{Jinlong Lei,  Peng Yi, Jie Chen, and Yiguang Hong
\thanks{The authors are with the  Department of Control Science and Engineering,
Tongji University,  Shanghai,   201804, China; }
\thanks{ Email address: leijinlong@tongji.edu.cn (J. Lei),  yipeng@amss.ac.cn (P. Yi),
 chenjie@bit.edu.cn(J. Chen) yghong@iss.ac.cn (Y. Hong).}
}

\begin{document}
\allowdisplaybreaks
\maketitle

%% ------------------------------------------------------------------
%% ABSTRACT
%% ------------------------------------------------------------------
\begin{abstract}
 The  paper considers    distributed  stochastic optimization  over randomly switching networks, where agents collaboratively minimize   the average of all agents' local  expectation-valued convex cost  functions. Due to the stochasticity in gradient observations, distributedness of local functions, and randomness of  communication topologies,
 distributed algorithms with a convergence guarantee under fixed step-sizes have not been achieved yet.
This work   incorporates variance reduction  scheme  into the distributed stochastic gradient tracking algorithm, where    local gradients are estimated by averaging  across a variable   number of  sampled gradients.
With an identically and independently distributed (i.i.d.) random network,
we   show that all agents' iterates  converge almost surely to the same optimal solution under fixed step-sizes. When the global  cost function is strongly convex and the sample size increases at a geometric rate, we  prove that the  iterates geometrically   converge   to the unique optimal solution, and establish the iteration, oracle, and  communication  complexity. The algorithm performance including rate and complexity analysis   are further  investigated with  constant step-sizes and  a polynomially increasing sample size. Finally, the empirical algorithm  performance  are illustrated with
numerical examples.
\end{abstract}

%\begin{IEEEkeywords}
%  proximal stochastic gradient,  randomized block-coordinate descent, nonsmooth optimization,  nonconvex stochastic optimization, variable sample-size schemes
%\end{IEEEkeywords}

\section{Introduction}

Distributed optimization has wide applications in  economic dispatch in power grids \cite{yi2016initialization,liu2019distributed}, trajectory planning and  control   for multi-robots \cite{fang2018optimization},   as well as machine learning over Internet of Things \cite{nedic2017fast,wang2018distributed,forero2010consensus}.
In distributed optimization, a group of agents connected over networks cooperatively minimizes the average of  all agents'  local cost functions. Prominent  first-order   distributed optimization algorithms have been developed, including    primal domain   methods that  combine classical (sub)gradient steps with  local averaging, such as  distributed subgradient methods \cite{nedic2009distributed,lou2017privacy}, first-order methods with historical gradients \cite{shi2015extra},
 distributed Nesterov gradient methods \cite{jakovetic2014fast,xin2018linear},  and distributed gradient tracking methods \cite{nedic2017achieving,xu2017convergence};  dual domain methods employing the Lagrangian dual,
  e.g., distributed dual decomposition \cite{palomar2006tutorial} and distributed ADMM \cite{mota2013d,shi2014linear};  and primal-dual domain methods \cite{chang2014distributed,yi2015distributed,lei2016primal}. \red{In addition, there are some works on second-order methods for stochastic and  distributed optimization \cite{schraudolph2007stochastic,mokhtari2015global,mokhtari2016dqm,wai2018sucag,li2020d}.}
Please refer to the survey   \cite{yang2019survey,notarstefano2019distributed} for the   recent progress.

Among various formulations in distributed optimization,
stochastic optimization has  particular research interests in multi-agent networks  due to its
applications in distributed estimation, stochastic control and  machine learning \cite{ram2010distributed, bianchi2012convergence,agarwal2011distributed}, where the local cost function is the expectation of a stochastic function.
In big data driven applications, the expectation is   a sum of sampling functions, while  it might be prohibitive or cumbersome to compute the exact gradient. Stochastic gradient descent (SGD)     becomes  popular since it is relatively easy to implement and scales well in large datasets \cite{bottou2018optimization}.

In distributed stochastic optimization,  each agent utilizes locally available sampled gradients and neighboring  information to cooperatively seek the optimal solution.
For nonsmooth convex cost functions, \cite{ram2010distributed} investigated a distributed stochastic subgradient
projection algorithm  and showed its mean convergence with  both the gradient estimation error and the step-size diminishing to zero, while \cite{srivastava2011distributed} further considered   the asynchronous distributed  SGD     over random networks  and   proved  almost sure convergence with  two diminishing step-size sequences.
For non-convex problems,  \cite{bianchi2012convergence}  showed that distributed SGD methods with diminishing step-sizes can guarantee the  almost sure  convergence to   Karush-Kuhn-Tucker points.
 Beyond the distributed SGD,  \cite{lei2018asymptotic} proposed  a  primal-dual method for distributed stochastic convex optimization over  random networks  corrupted with stochastic communication noises, and showed the  almost sure convergence with diminishing step-sizes.
{\it However, a convergent algorithm with non-diminishing step-sizes is
desirable in distributed stochastic optimization, since it can lead to a faster convergence rate, save the communication cost, and endow the multi-agent network with adaptivity under model drifting\cite{sayed2014adaptation}.}
As far as we know,   distributed algorithms for {\it stochastic convex optimization} with a convergence guarantee under non-diminishing step-sizes have not been  achieved yet.%, which is  one main motivation of this work.

There have been some distributed algorithms investigating   {\it strongly convex stochastic optimization}.  For example, \cite{nedic2016stochastic} proposed a subgradient-push method  over time-varying directed graphs
with  convergence rate   $\mathcal{O}(\ln(k)/k) $,
while  \cite{sayin2017stochastic} designed  a  stochastic subgradient descent with time-dependent averaging
and obtained a convergence  rate  $\mathcal{O}(1/k) $. In addition, \cite{yuan2018optimal}
considered a distributed stochastic mirror descent method with rate  $\mathcal{O}(\ln(k)/k) $
for  non-smooth functions. While for  random networks, \cite{jakovetic2018convergence} established  the mean-squared convergence rate  $\mathcal{O}(1/k) $ for distributed SGD.
Since the aforementioned  works \cite{jakovetic2018convergence,yuan2018optimal,sayin2017stochastic,nedic2016stochastic}
adopted diminishing step-sizes, the derived convergence rates are not
comparable with the  geometric   rate of  deterministic  strongly convex optimization with constant step-sizes.
Recently, \cite{pu2020distributed} and \cite{xin2019distributed}
proposed  distributed stochastic  gradient tracking methods
with constant step-sizes, but  only showed that the iterates  are attracted to a {\it neighborhood}
 of the optimal solution in expectation at a geometric rate.
With a different perspective,  \cite{alghunaim2019distributed} proposed a  distributed penalty gradient method for constrained stochastic optimization with a fixed step-size,  but also showed  the \red{geometric} convergence   to a {\it neighborhood} of the optimal solution. Thereby, how to distributedly achieve a linear convergence for strongly convex stochastic optimization needs further investigation.

 Variance reduction schemes  have gained increasing research interests in  stochastic convex
optimization~\cite{shanbhag2015budget,ghadimi2016accelerated,defazio2014saga,johnson2013accelerating}.
 In the class of variable sample-size schemes, the true gradient
is estimated by the  average of an increasing size of sampled gradients,
 which can progressively reduce the variance of the sample-averaged gradients. For example, \cite{shanbhag2015budget} obtained the geometric rate  for strongly convex problems, while \cite{ghadimi2016accelerated} combined the accelerated method  and  proved the rate  $\mathcal{O}(1/k^2)$ for  smooth convex  problems.
Alternative  variance reduction  schemes like SAGA \cite{defazio2014saga} and  SVRG \cite{johnson2013accelerating},   mainly applied to finite-sum optimization problems in machine learning,  lead  to the  recovery of the convergence rates
in deterministic cases. Such  schemes are also investigated in distributed finite-sum optimization \cite{mokhtari2016dsa,yuan2018variance,xin2019variance}, but {\it relying on periodically using   exact gradients}.
However, distributed variance reduced schemes for general distributed stochastic optimization without using exact gradients remains   open.

This paper aims to provide a fast and communication-efficient algorithm
for distributed stochastic optimization, where the communication   is many times  of local computation cost. We incorporate the   variable sample-size scheme   into  the distributed stochastic gradient tracking algorithm  \cite{pu2018push},
 and derive the following results.
\begin{itemize}
\item We propose a distributed   algorithm,
 where each agent estimates its local gradients  by a  variable number of sampled gradients, takes   a weighted averaging of   its neighbors'  iterates,  and moves towards   the
negative direction of the locally weighted combination of  its neighbours' gradient estimations.
%It is worth noting that the agents adopt   fixed but heterogeneous  step-sizes.
\item Assume that each sampled gradient   is   unbiased   with  a bounded variance, and each gradient   function  is Lipschitz continuous. For i.i.d.  random networks with connected mean graph,  we prove the almost sure convergence  for {\it merely convex} functions,
only requiring  the  sample size $N(k)$ satisfies  $\sum_{k=0}^{\infty} {1\over N(k)}<\infty$, which
is not necessary monotonically increasing.

\item If  the global   cost function is  strongly convex, we prove the   geometric   convergence with a geometrically increasing sample size, and  obtain  the iteration, communication, and oracle complexity $\mathcal{O}(\ln (1/\epsilon))$,  $\mathcal{O}(1/\epsilon^2)$, and $ \mathcal{O}(|\mathcal{E}|\ln(1/\epsilon))$
for  achieving  an $\epsilon$-optimal solution $ \mathbb{E}[||x-x^* ||]< \epsilon$.
%The iteration number $\mathcal{O}(\ln (1/\epsilon))$ matches that of the deterministic regimes,  the total number of  sampled gradients   $\mathcal{O}(1/\epsilon^2)$ is  of the  same order as that of centralized SGD, and the number of communication rounds is $ \mathcal{O}(|\mathcal{E}|\ln(1/\epsilon))$.
 We further show that with a constant sample size, the estimates
geometrically converge to a neighborhood of the optimal solution, and
investigate the polynomial rate   and  complexity bounds with a polynomially increasing sample size. The above results quantitatively
provide the trade-off between communication complexity and computation complexity  for distributed stochastic  optimization.
\end{itemize}

The novel perspective of this paper is  that  by  progressively  reducing {the   variance of  gradient   noises through increasing the sample size, we can adopt constant step-sizes to achieve an exact convergence in distributed stochastic optimization. Compared with algorithms with diminishing stepsizes in  \cite{ram2010distributed,srivastava2011distributed,bianchi2012convergence}, the proposed algorithm can achieve a faster convergence with constant step-sizes, hence  can significantly reduce the communication costs.
Moreover, for strongly convex stochastic optimization,
the derived iteration complexity is  of the same order as the  centralized algorithm in  {\it deterministic} cases \cite{nesterov2013introductory}.
The oracle complexity is also comparable with   centralized SGD, for example,
the bound of  \cite{hazan2014beyond} is $O(\frac{1}{\epsilon})$ for making the suboptimality gap $\mathbb{E}[F(x)]-F(x^*)<\epsilon$.
Compared with existing distributed stochastic optimization methods \cite{jakovetic2018convergence,yuan2018optimal,sayin2017stochastic,nedic2016stochastic,pu2020distributed,xin2019distributed} and \cite{mokhtari2016dsa,yuan2018variance,xin2019variance}, {\it the  proposed scheme saves the communication costs without increasing the overall sampling burden too much or using the exact gradient periodically.}

The paper is organized as follows.
A  distributed   variable sample-size stochastic gradient tracking algorithm  is proposed  in  Section II.
The almost sure convergence  for convex functions is provided in Section III.
Then the  geometric (resp. polynomial) convergence rate along with   complexity bounds are established in Section IV for strongly convex functions with geometrically (resp. polynomially) increasing sample size.
The numerical studies are presented  in Section V, while concluding remarks are given in Section VI.
In addition,  the proofs of   lemmas and theorems are presented in   Appendix.

 {\em Notations}.     Depending on the argument, $|\cdot|$ stands for the absolute value of a real number or the cardinality of a set.  The Euclidean norm of a vector or a matrix  is denoted as $\|\cdot\|_2$ or $\|\cdot\|$. The spectral radius of a matrix $A$ is denoted as $\red{\rho}(A)$.
Let $\otimes$ denote  the Kronecker  product.  The expectation  of a random variable is denoted as $\mathbb{E}[\cdot]$.
  Let  $ \mathbf{1}_n  $ denote   the  $n$-dimensional column vectors with all entries equal to 1 and $ I_d $  denote the $d\times d$ identity matrix.
  A directed graph  is denoted by  $\mathcal{G}=\{ \mathcal{V},\mathcal{E}\},$ where $\mathcal{V}=\{1,\dots,n\}$  is a finite set of  nodes and  an  edge  $(i,j) \in \mathcal{E}$ if node $j$ can receive information from agent $i$.
  A   directed path in $\mathcal{G}$  from $v_1$ to $v_{p }$ is a  sequence of distinct nodes, $v_1,  \dots, v_{p}$,
  such that  $(v_m, v_{m+1}) \in \mathcal{E}$ for all $m=1,\dots,p-1$. The graph $\mathcal{G}$ is termed {\it strongly    connected} if for any two distinct nodes $i,j\in\mathcal{V}$, there is a directed path from node $i$ to node $j.$
Given a nonnegative  matrix  $A=[a_{ij}] \in \mathbb{R}^{n \times n}$, denote by
    $\mathcal{G}_{A}= \{\mathcal{V},\mathcal{E}_A\}$  the corresponding digraph, where  $\mathcal{V }=\{ 1,\cdots,n\}$
       and   $(j,i)\in\mathcal{E}_A$ if   $a_{ij} >0 $.

\section{Problem  statement and distributed  algorithm}
In this section, we first formulate a distributed stochastic optimization problem.
Then we propose a  fully distributed  stochastic gradient tracking algorithm, which
used a  variable  number of   sampled gradients to estimate  exact gradients.

\subsection{Problem formulation}
We consider a network of   $n$ agents  indexed as $\mathcal{V}=\big\{1,\dots,n\big\}$.
Each agent $i \in \mathcal{V}$ has an expectation-valued  cost function  $f_i(x)\triangleq \mathbb{E}_{\xi_i}[h_i(x,\xi_i)]$, where  $x\in \mathbb{R}^d $,  the random vector  $\xi_i: { \Omega_i} \to \Real^{m_i}$ is
 defined on the probability space $({ \Omega_i}, {\cal F}_i, \mathbb{P})$,
and $h_i : \mathbb{R}^d\times \mathbb{R}^{m_i} \to \mathbb{R}$ is a proper scalar-valued  function. The    agents in the network need to cooperatively find an optimal solution that  minimizes the average of all agents' local cost functions, i.e.,
  \begin{equation}\label{problem1}
\begin{split}
\min_{x\in \mathbb{R}^d }  F(x )\triangleq {1\over n } \sum_{i=1}^n f_i(x).
\end{split}
\end{equation}
%We aim to design  a   distributed algorithm to drive all agents'  iterates to the optimal solution, explore  its rate of convergence, and establish  the  complexity bounds  for obtaining an optimal solution with a prescribed accuracy.

The discrete time is slotted at $k=0,1,2,\dots$. The interaction  among the   agents  at time $k$  is described by a directed     graph $\mathcal{G}(k)=\{ \mathcal{V},\mathcal{E}(k)\}$, where  $(j,i) \in \mathcal{E}(k)$ if node $i$ can receive information from agent $j$ at time $k$.  Denote by  $\mathcal{N}_i(k)\triangleq \{j\in \mathcal{V}: (j,i)\in \mathcal{E}(k)\}$  the neighboring  set of node $i$ at time $k$. The  corresponding  adjacency matrix   is   $A(k)=[a_{ij}(k)]_{i,j=1}^n$, where $a_{ij}(k)>0$  if $( j,i) \in \mathcal{E}(k)$   and   $a_{ij}(k)>0$ ,  otherwise.
Below are the assumptions on the communication graphs.
 \begin{assumption}\label{ass-graph}
 (i) At each time $k\geq 0,$    $A(k)$ is doubly stochastic, i.e.,
 $\mathbf{1}_n^T A(k)=\mathbf{1}_n^T$ and  $ A(k)\mathbf{1}_n=\mathbf{1}_n$.
% \begin{align}
%  \sum_{i=1}^n a_{ij}(k)=1,~  \sum_{j=1}^n a_{ij}(k)=1, \quad \forall j\in\mathcal{V},   \label{adj-a}   .
% \end{align}

 \noindent(ii) $\{A(k)\}_{k\geq 0}$ is an i.i.d. matrix  sequence.

\noindent (iii)   The   graph    $\mathcal{G}_{ \bar{A}}$ generated by the expected adjacency matrix
$ \bar{A}  = \mathbb{E}[A(k) ]$ is  strongly connected.

\end{assumption}

\blue{
\begin{remark}Random graphs modelled in Assumption \ref{ass-graph} can cover  i.i.d. undirected graphs \cite{Lobel2011},  random gossip and broadcast communications \cite{srivastava2011distributed,jakovetic2018convergence}, etc.   Assumption \ref{ass-graph}(i) requires each digraph  $\mathcal{G}(k)$ to be weight-balanced,
which is also used in existing works, such as \cite{xu2017convergence,xin2018linear,nedic2009distributed}.
Specifically, for the gossip scheme in undirected and connected  underlying graphs, the doubly stochastic adjacency matrix was  designed \cite{NedicOR18}. Nevertheless, it is usually non-trivial to generate doubly stochastic weights for general digraphs, but there are distributed algorithms to fulfill the task, e.g., \cite{CortesG12}.
Assumption \ref{ass-graph}(ii) implies that   the graph sequence $\{\mathcal{G}(k)\}$ is independent and identically distributed over time $k$. Assumption \ref{ass-graph}(iii) imposes a mild connectivity condition
among   agents that in expectation,  an agent can receive the information from every other agent directly or indirectly through a directed path.
\hfill $\Box$

\end{remark}}

We   require the cost functions to be convex and smooth.
\begin{assumption}\label{ass-convex}  For each agent $i\in \mathcal{V},$\\
(i) the cost  function $f_i$ is convex;\\
(ii) the gradient function    $\nabla f_i$ is  $L$-Lipschitz continuous, i.e.,
   $$  \| \nabla f_i(x_1)-  \nabla f_i(x_2) \|  \leq L \| x_1-x_2\|,\quad \forall x_1,x_2 \in \mathbb{R}^d.$$
\end{assumption}

Denote by $X^*\triangleq \{ x \in \mathbb{R}^d: \nabla F(x )=0\}$ the optimal solution set and by $F^*$
the optimal function  value. By the first-order optimality condition, the optimal solution $x^* \in X^*$ satisfies $\nabla F(x^*)=0.$

Suppose that for agent $i \in \mathcal{V}$, there exists a {\em stochastic first-order oracle}
that returns a  sampled gradient $\nabla_{x} h_i(x,\xi)$ given $x,\xi$, which is an {\it unbiased}  estimator of  $\nabla f_i(x)$ with a bounded second-order moment.
%Next follows the assumption on the stochastic first-order oracle.

  \begin{assumption}\label{ass-noise}  There exists a constant $\nu>0$ such that
  for each $i\in \mathcal{V}$  and any given $x\in \mathbb{R}^d$,
  $\mathbb{E}_{\xi_i}[\nabla_x h_i(x,\xi_i)]=\nabla f_i(x) $ and $
 \mathbb{E}_{\xi_i}[\| \nabla_x h_i(x,\xi_i)-\nabla f_i(x)\|^2 |x]  \leq \nu^2 .$
  %\begin{align*}
%&  \mathbb{E}_{\xi_i}[\nabla_x h_i(x,\xi_i)]=\nabla f_i(x),  {~\rm and ~} \\&
% \mathbb{E}_{\xi_i}[\| \nabla_x h_i(x,\xi_i)-\nabla f_i(x)\|^2 |x]  \leq \nu^2 .
%  \end{align*}
    \end{assumption}

\subsection{Distributed   algorithm with variable sample-sizes }

Each agent $i$ at time $k$ maintains two estimates $x_i(k)$  and   $y_i(k)$,
which are  used  to estimate the optimal solution and to  track  the average gradient, respectively.
 Since  the  exact gradient of each expectation-valued cost function $f_i(x)$ is unavailable, we  approximate it by  averaging through  a variable number of  sampled gradients,  \begin{align}\label{est-grad}
 \tilde{g}_i(x_i(k))={1\over N(k)} \sum_{p=1}^{N(k)}\nabla_{x}  h_i(x_i(k),\xi^p_i(k)), \quad \forall k\geq 0,
 \end{align}
where   $N(k)$   is the number of sampled gradients  utilized at  time   $k $ and    the  samples $\{ \xi^p_i(k)\}_{p=1}^{N(k)}$ are   randomly and independently  generated  from the  probability space $({ \Omega_i}, {\cal F}_i, \mathbb{P})$.
The gradient estimate  given by \eqref{est-grad} is an unbiased estimate of the exact gradient $\nabla f_i(x_i(k)) $,
and the  variance of the gradient observation noise $ \tilde{g}_i(x_i(k)) - \nabla f_i(x_i(k)) $ can be progressively reduced by increasing  the sample size $N(k)$. By  combining the distributed gradient tracking scheme     with such a variance reduction scheme, we  obtain  Algorithm~\ref{alg-sto-gradient_track}.
%We will specify the selection of the constant stepsizes  $\alpha_i,i\in \mathcal{V}$
% and the sample size $N(k)$ upon convergence analysis.
%

\begin{algorithm}
\caption{Distributed variable sample-size stochastic gradient tracking algorithm}\label{alg-sto-gradient_track}
{\em Initialization}: Set $k:=0$. For any $i = 1, \hdots, n$, let  $y_i(0)= \tilde{g}_i(x_i(0))$ with  arbitrary initial  $ x_i(0) \in  \mathbb{R}^{d} $.

{\em Iterate until convergence.}

Each agent  $i=1,\cdots, n$ updates  its estimates as follows,
\begin{subequations}
\begin{align}
 x_i(k+1)& =\sum_{j\in \mathcal{N}_i(k)}   a_{ij}(k)  x_j(k)-\alpha_i y_i(k), \label{push}\\
 y_i(k+1)&= \sum_{j\in \mathcal{N}_i(k)}   a_{ij}(k)  y_j(k) +\tilde{g}_i(x_i(k+1))-  \tilde{g}_i(x_i(k)), \label{pull}
  \end{align}
\end{subequations}
 where   $\alpha_i>0$ is the {\it fixed} step-size used by agent $i$,
  and  $ \tilde{g}_i(x_i(k))$ is given in \eqref{est-grad}.
   \end{algorithm}

Note that for each agent $i\in \mathcal{V},$ the  implementation of      \eqref{push}   requires  its neighbors'  estimates of the optimal solution  $\{ x_j(k)\}_{ j\in \mathcal{N}_i(k)}$,
 while the update of $y_i(k+1) $    characterized by     \eqref{pull}  uses its  local gradient estimate  as well as  its neighbors'  information $\{ y_j(k)\}_{ j\in \mathcal{N}_i(k)}$ to  asymptotically track   the dynamical average  gradient across the network.
  Therefore,  Algorithm~\ref{alg-sto-gradient_track} is fully distributed  since  each agent merely relies on its local samples and  its neighboring agents'   information.

\section{Almost Sure Convergence for Convex Functions }\label{sec:as}

In this section, we provide the almost sure convergence  of the algorithm for merely convex cost functions.
\subsection{Preliminary lemmas}
Define  the  gradient observation noise  as follows,
\begin{equation}\label{def-w}
\begin{split}
  &w_i(k)\triangleq  \tilde{g}_i(x_i(k)) -\nabla f_i(x_i(k)),
  \\  & w(k ) \triangleq \big(w_1(k )^T, \cdots, w_n(k)^T \big)^T \in \mathbb{R}^{nd}.
\end{split}\end{equation}
Denote by
\begin{equation}\label{def-fw}
\begin{split}
  &x (k) \triangleq \big( x_1(k) ^T, \cdots, x_n(k) ^T \big)^T \in \mathbb{R}^{nd},
    \\& y(k ) \triangleq \big(y_1(k )^T, \cdots, y_n(k)^T \big)^T   ,
      \\ & \red{\nabla f(k) }\triangleq \big(\nabla f_1(x_1(k)) ^T, \cdots, \nabla f_n(x_n(k)) ^T \big)^T  ,
\\ & {\rm~and~}\bm{\alpha}\triangleq {\rm diag} \{\alpha_1,\cdots, \alpha_n\}\in \mathbb{R}^{n \times n} ,
\end{split}
\end{equation}
where ${\rm diag} \{\alpha_1,\cdots, \alpha_n\}$ denotes a diagonal matrix with $\alpha_i$ in the $i$th diagonal.
Then  Algorithm \ref{alg-sto-gradient_track} can be written in a compact form as follows,
 \begin{subequations}
\begin{align}
 x (k+1)& =(A(k)  \otimes I_d ) x(k)- (\bm{\alpha}  \otimes I_d)  y(k) , \label{push2}\\
 y(k+1)&= (A(k)  \otimes I_d )   y(k) +\red{\nabla f(k+1) }  \notag\\&\quad +w(k+1)-\red{\nabla f(k) }-w(k ).\label{pull2}
  \end{align}
\end{subequations}

Denote the averaged estimate of the optimal  solution and  the averaged gradient across the network  as
\begin{align}\label{def-xy}
\bar{x}(k)={1\over n} \sum_{i=1}^n x_i(k)  {\rm~and}~\bar{y}(k)={1\over n} \sum_{i=1}^n y_i(k) .
\end{align}
We further denote  by $D_{\bot}=I_n-{\mathbf{1}_n\mathbf{1}_n^T \over n},$ and by
\begin{align}
 & \tilde{x}(k) \triangleq (D_{\bot}\otimes I_d)  x(k)= x(k)- (\mathbf{1}_n \otimes I_d)  \bar{x}(k) ,\label{def-tidlex}
 \\&  \tilde{y}(k)\triangleq  (D_{\bot}\otimes I_d)  y(k)=y(k)- (\mathbf{1}_n \otimes I_d)  \bar{y}(k).\label{def-tidley}
\end{align}
Define $\mathcal{F}(k)\triangleq \big \{x(0),A(0),\cdots, A(k-1), \{\xi^p_i(t)\}_{p=1}^{N(t)}, 0\leq t\leq k,i=1,\cdots, n \big \}.$ From Algorithm \ref{alg-sto-gradient_track}  it is seen that both $x(k)$ and $y(k)$ are adapted to  $\mathcal{F}(k)$, hence    $\tilde{x}(k)$ and $\tilde{y}(k)$ are adapted to  $\mathcal{F}(k)$.

Recall that $ \tilde{x}(k) $ is adapted to  $\mathcal{F}(k)$ and $A(k)$ is independent of $\mathcal{F}(k).$ Then by using $A(k){\mathbf{1}_n\mathbf{1}_n^T \over n}={\mathbf{1}_n\mathbf{1}_n^T \over n} $ and  $A(k)^T{\mathbf{1}_n\mathbf{1}_n^T \over n}={\mathbf{1}_n\mathbf{1}_n^T \over n},$  we derive
\begin{align*}
&\mathbb{E}\big[\| (A(k)-\mathbf{1}_n\mathbf{1}_n^T / n) \otimes I_d\tilde{x}(k)\|^2 | \mathcal{F}(k) \big]
\\&=\mathbb{E}\big[\tilde{x}(k)^T (A(k)-\mathbf{1}_n\mathbf{1}_n^T / n)^T (A(k)-\mathbf{1}_n\mathbf{1}_n^T / n)\otimes I_d \tilde{x}(k)  | \mathcal{F}(k) \big]
\\&=\tilde{x}(k)^T  \left(\mathbb{E} [  A(k)^TA(k)] -{\mathbf{1}_n\mathbf{1}_n^T \over n}  \right) \otimes I_d \tilde{x}(k)
\\& \leq \rho \left(\mathbb{E} [  A(k)^TA(k)] -{\mathbf{1}_n\mathbf{1}_n^T \over n}  \right)
\|\tilde{x}(k)\|^2.
\end{align*}
Therefore,  by applying  the   Jensen's inequality for  conditional expectations, we obtain
that 
\begin{equation}\label{bd-sr}
\begin{split}
&  \mathbb{E}\big[\| (A(k)-\mathbf{1}_n\mathbf{1}_n^T / n) \otimes I_d\tilde{x}(k)\| | \mathcal{F}(k) \big]
\\&\leq \rho_1 \| \tilde{x}(k)\| {\rm~with~} \rho_1\triangleq \sqrt{\rho (\mathbb{E} [  A(k)^TA(k) -{\mathbf{1}_n\mathbf{1}_n^T \over n} ] )}.
\end{split}
\end{equation}

By Assumption \ref{ass-graph}, we see that the graph generated by the matrix $\mathbb{E} [  A(k)^TA(k) ]$ is  undirected and  connected. Thus, $ \rho (\mathbb{E} [  A(k)^TA(k) -{\mathbf{1}_n\mathbf{1}_n^T \over n} ] )\in (0,1)$.
The parameter $ \rho_1$ depends  on the network topology, where larger $\rho_1$ implies worse network connectivity.
 It was shown in \cite[Proposition 5]{NedicOR18} that when the  weight of adjacency matrix following the Lazy Metropolis rule,
  $1-\rho_1=\mathcal{O}(1/n^2)$ for path or star graph, $1-\rho_1=\mathcal{O}(n^{-1})$ for lattice graph,
  etc. For general classes of graphs, there are some distributed algorithms to estimate the network connectivity
for fixed graphs, such as \cite{FranceschelliGGS13}.

Define
\begin{equation}\label{def-alpha}
\begin{split}
&\bar{\alpha} \triangleq {1\over n} \sum_{i=1}^n \alpha_i,~ c_1\triangleq \sqrt{\sum_{i=1}^n (\alpha_i-\bar{\alpha})^2}, \\
& \alpha_{\max} \triangleq \max_{i\in \mathcal{V}} \alpha_i,~ {\rm and} \; c_2\triangleq \sqrt{\sum_{i=1}^n \alpha_i^2}.
 \end{split}
 \end{equation}
In the following lemma,  We obtain  the following  bounds on the consensus errors  $\tilde{x}(k)$ and $\tilde{y}(k)$.

\begin{lemma}\label{lem-cond} Suppose Assumptions \ref{ass-graph} and \ref{ass-convex} hold.
With  Algorithm  \ref{alg-sto-gradient_track},    we have that, for any $k\geq 0 $,
\begin{align}
 &\mathbb{E}[\| \tilde{x} (k+1)  \| | \mathcal{F}(k)] \notag
 \\&\leq  \rho_1 \| \tilde{x}(k)\|
  + c_1 \| \bar{y}(k)\| + \alpha_{\max} \| \tilde{y}(k)\|,\label{recur-tildex}
  \end{align}
  and
  \begin{align}
&\mathbb{E}[\| \tilde{y} (k+1)  \| | \mathcal{F}(k)] \leq     ( \rho_1 +\alpha_{\max}L) \| \tilde{y}(k)\| +   L \| \tilde{x}(k) \|\notag
\\&     +c_2L \|    \bar{y}(k)  \| + \mathbb{E}[ \| w(k+1)-w(k ) \| | \mathcal{F}(k)].\label{recur-yerr2}
\end{align}
%\red{where $\rho_1$ is defined in \eqref{bd-sr}.}
\end{lemma}
\begin{IEEEproof}  The proof is given  in Appendix \ref{app:A1}. \end{IEEEproof}

\vskip 3mm
\begin{lemma}\label{lem-exp}   Suppose Assumptions \ref{ass-graph},  \ref{ass-convex}, and \ref{ass-noise} hold.
Consider Algorithm  \ref{alg-sto-gradient_track}, where  $\alpha_i< { (1-\rho_1)^2\over (2-\rho_1)L}$ for each $i\in \mathcal{V}$. Define $e(k) \triangleq \sqrt{\mathbb{E}  [\| \tilde{x} (k )  \|^2+ \| \tilde{y} (k )  \|^2] }$.
 Then $$\rho_2\triangleq {2\rho_1+\alpha_{\max}L + \sqrt{\alpha_{\max}^2L^2+4\alpha_{\max}L} \over 2}<1,$$
 and  the following holds with $p_k \triangleq  { \sqrt{n} \nu  \over \sqrt{N(k+1)}} +{ \sqrt{n} \nu  \over \sqrt{N(k)}} :$
\begin{align}
\sqrt{\sum_{s=0}^K  e(s)^2 }& \leq \sqrt{3   \over 1-\rho_2^2} e(0) + {  \sqrt{3  }  \over 1-\rho_2}  \sqrt{\sum_{s=0}^{K } p_s^2}\notag
\\&+  { \sqrt{3( c_1^2+   c_2^2L^2)} \over  1-\rho_2} \sqrt{\sum_{s=0}^{K }\mathbb{E}[\| \bar{y}(s)\|] ^2} .
\end{align}
\end{lemma}

\begin{IEEEproof}
  The proof is given   in Appendix \ref{app:A2}. \end{IEEEproof}

\subsection{Almost sure convergence}
Next, we give the almost sure convergence of Algorithm \ref{alg-sto-gradient_track}.

\begin{theorem}\label{thm-as}
 Suppose Assumptions \ref{ass-graph},  \ref{ass-convex}, and \ref{ass-noise} hold.
Let $\{x(k)\}$ and $\{y(k)\}$ be generated by Algorithm  \ref{alg-sto-gradient_track}, where
 $\sum_{k=0}^{\infty} {1\over N(k)}<\infty$. Then there exist sufficiently small  $\alpha_i>0,i\in \mathcal{V}$, which
 possibly \red{depends} on  $\rho_1$, $L,$ and
 \begin{equation}\label{def-da} d_{\alpha}\triangleq {\sqrt{\sum_{i=1}^n (\alpha_i-\bar{\alpha})^2} \over \sqrt{n} \bar{\alpha}}, \end{equation}
  such that
 \begin{equation}\label{thm1-result}
 \begin{split}
 &\lim_{k \rightarrow \infty  } \| \bar{x}(k )-x_i(k)\|=0, ~\forall i\in \mathcal{V}, \quad a.s.,
 \\&\lim_{k \rightarrow \infty  }F(\bar{x}(k ))=F^*  ,\quad a.s.
 \end{split}
 \end{equation}
\end{theorem}

\begin{IEEEproof}
  The proof is given   in Appendix \ref{app:A3}. \end{IEEEproof}

\begin{remark}Theorem \ref{thm-as} shows that the exact convergence in an almost sure  sense can be achieved for convex problems with constant step-sizes by adaptively choosing the sample size.
The  proposed algorithm  with constant step-sizes can achieve a faster convergence rate compared with the algorithms with diminishing step-sizes   \cite{jakovetic2018convergence,yuan2018optimal,sayin2017stochastic,nedic2016stochastic}.
As similar discussions in \cite{xu2015augmented}, one major reason
for considering the agent-specific stepsize is due to the heterogeneity
of agents and lacking of coordination involved in distributed computation. Theorem \ref{thm-as}  validates
that the distributed variable sample-size stochastic gradient tracking
algorithm with uncoordinated constant stepsizes can also achieve the
exact convergence to an optimal solution in the almost sure sense.

  Theorem \ref{thm-as}   uses the same sample size  just for the ease of   proof presentation. Suppose   agents utilize different sample size, i.e., agent $i$ uses $N_i(k)$ at time $k$.
Denote  by $N_{\min}(k)=\min\{N_i(k), i\in\mathcal{ V}\}$. Then  the  condition  $\sum_{k=0}^{\infty} {1\over N(k)}<\infty$ can be replaced with $\sum_{k=0}^{\infty} {1\over N_{\min}(k)}<\infty$.
 There are many ways for choosing  the batch-size $N(k)$, for example    $   k \ln^2(k)  $
 or $   k^{1+\delta}  $ with $\delta>0$. \hfill $\Box$
\end{remark}

The following corollary  gives a sufficient condition on  constant step-sizes
when all   agents take an identical  step-size. It quantitatively characterizes  the dependence on the  Lipschitz constant $L$ and the network connectivity parameter $\rho_1$.
 It can be seen that a larger Lipschitz constant  $L$ leads to a smaller upper bound
 of the  step-size, while a better network connectivity (i.e., smaller $\rho_1$) implies a larger   step-size.

\begin{corollary}\label{cor1}  Suppose Assumptions \ref{ass-graph},  \ref{ass-convex}, and \ref{ass-noise} hold.
Consider  Algorithm  \ref{alg-sto-gradient_track} with $\alpha_i\equiv \alpha $ and
 \begin{align}\label{def-alpha1} \alpha \in \left(0, \tfrac{ c_0+1+2\sqrt{3}L-\sqrt{(c_0+1+2\sqrt{3}L)^2-4c_0}}{ 2L} \right)\end{align} with $c_0\triangleq { (1-\rho_1)^2\over (2-\rho_1) }$. Then the results  established in     \eqref{thm1-result} hold.
\end{corollary}
 \begin{IEEEproof}
  The proof is given   in Appendix \ref{app:A4}. \end{IEEEproof}

%
%\begin{remark} For ease of    stepsizes selection, we  let all the agents take the same stepsize, i.e.,    $\alpha_i\equiv  \alpha $. Then $c_1=0$  and
%$q_2={\sqrt{3(c_1^2+c_2^2L^2)} \over 1-\rho_2}={\sqrt{3n }  \alpha L  \over  1-\rho_2} $. Thus,  from the proof in Appendix \ref{app:A3} we see that   $ 1  -    \alpha  L -   {\sqrt{3   }\alpha L^2 \over  1-\rho_2}  >0$
% is  a sufficient condition for the almost sure  convergence in Theorem  \ref{thm-as}.
%This combined  with the definition of $\rho_2$ in Lemma \ref{lem-cond} indicates that
%  $\alpha\in \left(0, { c_0+1+2\sqrt{3}L-\sqrt{(c_0+1+2\sqrt{3}L)^2-4c_0} \over 2L} \right)$ is a sufficient condition, where $c_0\triangleq { (1-\rho_1)^2\over (2-\rho_1) }$. Such a  parameter selection quantitatively characterizes  the influence of the Lipschitz constant $L$ and the network connectivity parameter $\rho_1$.   \hfill $\Box$
%\end{remark}

%Define $\beta=\alpha L < { (1-\rho_1)^2\over (2-\rho_1) }.$ Then $\rho_2=\rho_1+{\beta+\sqrt{\beta^2+4\beta} \over 2}$. Therefore,
%\begin{align*}
%&(1  -    \beta)\left(  { (1-\rho_1)^2\over (2-\rho_1) }- \beta     \right) >    2\sqrt{3  }\beta L \\
%  &\Rightarrow (1  -    \beta)\left(  1-\rho_1-{\beta+\sqrt{\beta^2+4\beta} \over 2} \right) >    \sqrt{3  }\beta L
%  \\& \Leftrightarrow     (1  -    \beta)(  1-\rho_2 ) -    \sqrt{3  }\beta L     >0
%\end{align*}

\section{Rate Analysis for Strongly Convex Functions}\label{sec:sc}

This section  explores the    convergence properties of Algorithm~\ref{alg-sto-gradient_track}  when the
 global cost function is strongly convex.
The geometric (resp. polynomial) convergence rate  is obtained if the number
of the sampled gradients  increases  at  a geometric (resp. polynomial)  rate. In addition,
  the complexity bounds for obtaining an $\epsilon$-optimal solution are established as well.

\subsection{Linear convergence rate analysis}

 \red{\begin{assumption}\label{ass-fun}
    The global cost function $F(x)$ is $\eta$-strongly convex, i.e.,  for any $x_1,x_2 \in \mathbb{R}^d$,
 $$(\nabla F(x_1)-  \nabla F(x_2) )^T(x_1-x_2) \geq \eta \| x_1-x_2\|^2.$$   \end{assumption}}

With Assumption \ref{ass-fun}, the problem \eqref{problem1} has a unique optimal solution, denoted by  $x^*$, and  $\nabla F(x^*)=0$. We analyze  the algorithm performance  by   characterizing the   interactions  among three   error sequences: (i) the  distance from the average estimate to   the optimal solution $\| \bar{x}(k) -x^*\|;$
(ii) the consensus error $  \| x(k) - (\mathbf{1}_n \otimes I_d)\bar{x}(k)  \| $;
and  (iii) the consensus  error  of the gradient trackers   $\| y(k) - (\mathbf{1}_n \otimes I_d)  \bar{y}(k)\|$.
We   will   bound the three  error sequences in terms of linear combinations of their past values in the following lemma.

\begin{lemma}\label{lem1}  Suppose Assumptions \ref{ass-graph}, \ref{ass-convex}(ii), and \ref{ass-fun} hold.   Consider  Algorithm \ref{alg-sto-gradient_track} with $0<\alpha_i \leq {2\over \eta+L}.$ Define
\begin{equation} \label{def-z}
\begin{split} z(k) & \triangleq
\begin{pmatrix}
\mathbb{E}[\| \bar{x}(k) -  x^*  \|] \\
\mathbb{E}[\| x(k) - (\mathbf{1}_n \otimes I_d)\bar{x}(k)\|]  \\
\mathbb{E}[ \| y(k) - (\mathbf{1}_n \otimes I_d)  \bar{y}(k)\|]
\end{pmatrix} ,  {\rm~and~} \\
J(\bm{\alpha}) & \triangleq
 \begin{pmatrix}
1-\bar{\alpha }\eta  & {\bar{ \alpha}L \over \sqrt{n}}  & {c_1 \over n}
 \\ c_1 L&  \rho_1+ {c_1L \over \sqrt{n} } & \alpha_{\max}  \\
 c_2 L^2 &L+{ c_2L^2 \over \sqrt{n}} &\rho_1 +\alpha_{\max}L \end{pmatrix} ,
 \end{split}
 \end{equation}
 where $\bar{\alpha} , c_1 ,\alpha_{\max}, c_2  $ are defined in \eqref{def-alpha}.
 Then   the following component-wise linear  matrix  inequality  holds for  any $k\geq 0$,
\begin{equation}\label{recursion-z0}
\begin{split}
& z(k+1) \leq  J(\bm{\alpha}) z(k)    \\& +\begin{pmatrix}
 {\bar{\alpha }\over \sqrt{n}}  \mathbb{E}[ \|   w(k)   \|] \|\\
 {c_1\over \sqrt{n}}  \mathbb{E}[ \|   w(k)   \|] \|\\ \mathbb{E}[\| w(k+1)-w(k ) \|] + {c_2 L \over \sqrt{n}} \mathbb{E}[ \|   w(k)   \|]
\end{pmatrix} .
\end{split}
\end{equation}
\end{lemma}
\begin{IEEEproof} The proof  can be found  in Appendix \ref{app:B1}. \end{IEEEproof}

Next, we show the geometric   convergence of Algorithm \ref{alg-sto-gradient_track}
 with geometrically increasing sample size and suitably selected step-sizes.
 \red{For non-identical step-sizes $\alpha_i$, we have   $d_{\alpha}>0$ by  the definition \eqref{def-da}.
Define $\kappa\triangleq L/\eta$ and let $\alpha_i$ satisfy the following with $\rho_1$ defined by \eqref{bd-sr}:
 \begin{align}\label{bd-stepsize}
 &0<\alpha_iL < \min \left\{ \beta^*,  {1- \rho_1 \over d_{\alpha} \kappa(L+\eta) } \right\}, \quad \forall i\in \mathcal{V} ,\\
 & {\rm~where~} \beta^* \triangleq {c_4+ \sqrt{c_4^2+4c_3(1- \rho_1)^2 }\over 2c_3}  \notag
 \end{align}
 whith $c_3= (\sqrt{d_{\alpha}^2+1}  )  (1 + \kappa d_{\alpha}^2    )+\kappa\sqrt{d_{\alpha}^2+1},$
 and $c_4=(1 +  (\kappa+1) d_{\alpha}  )(  1- \rho_1 )+1 + \kappa d_{\alpha}^2 +\kappa\sqrt{d_{\alpha}^2+1}    d_{\alpha}   (1- \rho_1) .$}

\begin{theorem}  \label{thm1}  Suppose Assumptions \ref{ass-graph}, \ref{ass-convex}(ii), \ref{ass-noise}, and \ref{ass-fun}    hold. Let  $\{x(k)\}$ and $\{ y(k)\}$ be generated by  Algorithm \ref{alg-sto-gradient_track} with $N(k)=\lceil q^{-2k} \rceil$ for some $q\in (0,1). $
  \red{Suppose  the  step-size   $\alpha_i ,i\in \mathcal{V}$  satisfies \eqref{bd-stepsize}, then  the spectral radius of  $   J(\bm{\alpha})$ in \eqref{def-z}, denoted by $\rho(J(\bm{\alpha})) $,  is strictly smaller than 1.}
  In addition, the  error sequence $z(k)$   converges    to zero at a linear  rate,  $$\mathcal{O}\big( \max\{ \rho(J(\bm{\alpha})) ,q \} ^k\big).$$
\end{theorem}

\begin{IEEEproof} The proof   can be found in Appendix \ref{app:B2}. \end{IEEEproof}
 \red{Eqn. \eqref{bd-stepsize} gives a sufficient  condition   for selecting    step-sizes  $\alpha_i ,i\in \mathcal{V}$  to guarantee that $\rho(J(\bm{\alpha})) <1$. It     shows how parameters $\eta, L,d_{\alpha}$ and $\rho_1$ influence the selection of constant  step-sizes.}
Theorem \ref{thm1}  implies that if  the number of sampled gradients is increased at a geometric rate $\lceil q^{-2k} \rceil$ with $q\in (0,1)$, the    error sequences   $ \mathbb{E}[\| \bar{x}(k) -x^*\| ] $ and $ \mathbb{E}[ \| x(k) - (\mathbf{1}_n \otimes I_d)\bar{x}(k)  \| ]$    converge  to zero  at a geometric rate.
\red{We   omitted the  big $O$ constant  in the statement of Theorem 2 due to its complicated expression. However, it is noticed from \eqref{def-z} that $J(\alpha)$   is just a $3\times 3$ matrix, which can be computed if the problem-related constants are given.
In this case, the explicit convergence rate can be computed with the inequality \eqref{recur-z} in Appendix B.B.}

The following corollary  \red{shows the   convergence rate for the case with an identical  step-size, i.e.,
 $\alpha_i\equiv\alpha$. Define
 \begin{align}\label{def-hatJ}
 \hat{J}( \alpha )  \triangleq
 \begin{pmatrix}
1- \alpha  \eta  & {  \alpha L \over \sqrt{n}}  & 0
 \\ 0&  \rho_1  & \alpha \\
 \sqrt{n}\alpha L^2 &L+{  \alpha L^2  } &\rho_1 +\alpha L \end{pmatrix} .
 \end{align}
The condition \eqref{stepsize}, making $\rho(\hat{J}(\alpha))<1$,   implies that
a better network connectivity (namely a smaller $\rho_1$)   leads to a larger $\alpha,$ while
the ill-conditioned optimization problem with a large $\kappa$   narrows the possible selection of $\alpha.$ }

\begin{corollary}\label{cor2}
Suppose Assumptions \ref{ass-graph}, \ref{ass-convex}(ii), \ref{ass-noise}, and \ref{ass-fun}    hold.
 Let  $\{x(k)\}$ and $\{ y(k)\}$ be generated by  Algorithm \ref{alg-sto-gradient_track},
 where \begin{align}\label{stepsize} \alpha_i\equiv\alpha
 < { 2-\rho_1 +\sqrt{(2-\rho_1)^2+4(1+\kappa)(1- \rho_1 )^2}\over 2L(1+\kappa)} . \end{align}
  \red{Set $N(k)=\lceil q^{-2k} \rceil$ for some $q\in (\rho(\hat{J}(\alpha)),1). $
 Then    \begin{equation*}
\begin{split}
&z(k)  \approx  \hat{J}(\alpha) ^k z(0) \\& +\nu q^{k-1} \big(I_3-  \hat{J}(\alpha)/q \big)^{-1} \begin{pmatrix}
   \alpha  \\ 0 \\  (1+q+  \alpha    L ) \sqrt{n}
\end{pmatrix}   .
\end{split}
\end{equation*} }
\end{corollary}
\begin{IEEEproof} The proof  is given in Appendix \ref{app:B3}. \end{IEEEproof}

 For strongly convex stochastic optimization, \cite{pu2020distributed,xin2019distributed,alghunaim2019distributed} also proved  geometric convergence rates but only to a neighborhood of the optimal solution. By progressively reducing the gradient noises with geometrically
  increasing batch-sizes, we  prove that   the exact and geometric convergence in a mean-squared sense.
  \blue{The following Corollary shows that when a constat sample size is used in Algorithm \ref{alg-sto-gradient_track},  the linear convergence to a neighborhood of the optimal solution
  can  be obtained as well. It can be seen that  the bounds  depend  on the network structure, batch-size and step-size, as well as the problem parameters $\eta, L, \nu.$ }

\blue{\begin{corollary}\label{cor3}  Let Assumptions \ref{ass-graph}, \ref{ass-convex}(ii), \ref{ass-noise}, and \ref{ass-fun}    hold. Consider Algorithm \ref{alg-sto-gradient_track} with $N(k)\equiv B$ for some positive integer $B,$
   where     $\alpha_i \equiv \alpha,i\in \mathcal{V}$  satisfies \eqref{stepsize}. Then
 $ \sup_{l\geq k}\mathbb{E}[\| \bar{x}(l) -x^*\| ] $ and $\sup_{l\geq k}  \mathbb{E}[ \| x(l) - (\mathbf{1}_n \otimes I_d)\bar{x}(l)  \| ]$ converge to  $\limsup_{k\to \infty} \mathbb{E}[\| \bar{x}(k) -x^*\| ] $ and  $\limsup_{k\to \infty} \mathbb{E}[ \| x(k) - (\mathbf{1}_n \otimes I_d)\bar{x}(k)  \| ]$
 with    a geometric  rate $\mathcal{O}\big( \rho(\hat{J}(\alpha)) ^k\big)$.  Furthermore,
\begin{align*}
&\limsup_{k\to \infty} \mathbb{E}[\| \bar{x}(k) -x^*\| ]
\\&\leq   {   \nu ( (1-\rho_1)^2+\rho_1\alpha L) \over \sqrt{B} \eta \left( (1- \rho_1 )^2-(1+\kappa) \alpha^2L^2 -  (2-\rho_1)\alpha L \right)}
 \end{align*} and
 \begin{align*}& \limsup_{k\to \infty} \mathbb{E}[ \| x(k) - (\mathbf{1}_n \otimes I_d)\bar{x}(k)  \| ]
\\& \leq{\alpha \sqrt{n} \nu  (\alpha L^2  +  \eta (2 + \alpha   L  )     ) \over \sqrt{B}\eta \left( (1- \rho_1 )^2-(1+\kappa) \alpha^2L^2 -  (2-\rho_1)\alpha L \right)}.
 \end{align*}
\end{corollary}
\begin{IEEEproof} The proof  is given in  Appendix \ref{app:B4}. \end{IEEEproof}}
\subsection{Complexity analysis}

Based on the geometric convergence rate established in Theorem \ref{thm1},
we are able  to  establish the complexity bounds
for obtaining an  $\epsilon$-optimal solution satisfying   $  \| z(k)\|   \leq \epsilon.$
 The   iteration  complexity is defined as  $K(\epsilon)$  such that $  \| z(k )\|   \leq \epsilon $ for any $ k \geq K(\epsilon)$. The oracle complexity,  measured by the total number of sampled gradients
 for deriving an   $\epsilon$-optimal solution,   can be computed as   $ \sum_{k=0}^{ K(\epsilon)} N(k).$
%The following theorem gives the complexity bounds.

\begin{theorem}  \label{thm2}  Let Assumptions \ref{ass-graph}, \ref{ass-convex}(ii), \ref{ass-noise}, and \ref{ass-fun} hold.
Consider  Algorithm \ref{alg-sto-gradient_track} with    $N(k)=\lceil q^{-2k} \rceil$ for some $q\in (0,1), $
 \red{where  the  step-size   $\alpha_i ,i\in \mathcal{V}$  satisfies \eqref{bd-stepsize}.}
\\ (i) When  $ \rho(J(\bm{\alpha}))<q<1,$ the
iteration  and oracle complexity  required  to obtain  an $\epsilon$-optimal solution  are $ \mathcal{O}(\ln(1/\epsilon))$ and $\mathcal{O}(1/\epsilon^2)$, respectively.
 \\ (ii) When $ 0<q<\rho(J(\bm{\alpha}))  $, the iteration and oracle complexity  required  to obtain  an $\epsilon$-optimal solution  are $ \mathcal{O}(\ln(1/\epsilon))$ and  $   (1/ \epsilon)^{ 2\ln(1/q) \over \ln(1/\rho(J(\bm{\alpha}))) }$, respectively.
\end{theorem}
\begin{IEEEproof}    (i).  $\rho(J(\bm{\alpha}))<q$.
With Theorem \ref{thm1}, there exists $C_1>0$ such that $\| z(k )\| \leq C_1q^k.$
Then for any $  k\geq  K_1(\epsilon) =  \ln(C_1/\epsilon) {1\over \ln(1/q)}$, we have $\| z(k )\| \leq   \epsilon .$
 This allows us to bound  the oracle complexity   by
\begin{align*}
&\sum_{k=0}^{K_1(\epsilon)} N(k)=\sum_{k=0}^{K_1(\epsilon)}  q^{-2k} \leq { q^{-2(K_1(\epsilon)+1)} \over q^{-2}-1}
\\& \leq  {1\over 1-q^2} q^{-2 { 1\over \ln(1/q) }   \ln(C_1/\epsilon) }
=  {1\over 1-q^2}  e^{\ln(q^{-2 }){ 1\over \ln(1/q) }   \ln(C_1/\epsilon) }
 \\&=  {1\over 1-q^2}  e^{2  \ln(C_1/\epsilon) }={C_1^2\over (1-q^2)\epsilon^2}.
\end{align*}

\noindent   (ii).    $\rho(J(\bm{\alpha}))>q$.
  With Theorem \ref{thm1}, there exists $C_2>0$ such that $\| z(k )\| \leq C_2\rho(J(\bm{\alpha}))^k.$
   Then for any $k\geq K_2(\epsilon) \triangleq { 1\over \ln(1/\rho(J(\bm{\alpha}))) }  \ln  \big( { C_2 \over \epsilon}\big)  $,  we have
   $\| z(k )\| \leq   \epsilon .$ This allows us to bound the  oracle complexity   by
\begin{align*}
&\sum_{k=0}^{K_2(\epsilon)} N(k) \leq  {   q^{-2 (K_2(\epsilon)+1)} \over q^{-2}-1}
\leq   {1\over 1-q^2}  q^{-2 { 1\over \ln(1/\rho(J(\bm{\alpha}))) }  \ln  \left( { C_2 \over \epsilon}\right) }
\\&=    {e^{\ln(q^{-2 }){ 1\over \ln(1/\rho(J(\bm{\alpha}))) }  \ln  \left( { C_2 / \epsilon}\right) }\over 1-q^2}
 =  {1\over 1-q^2}    \left( {C_2 \over \epsilon}\right)^{2\ln(1/q) \over \ln(1/\rho(J(\bm{\alpha}))) }.
\end{align*}
\end{IEEEproof}

\begin{remark}\label{remark4}
Theorem \ref{thm2}  shows that for   geometrically increasing batch-size, the
number of iterations required to obtain an $\epsilon$-optimal solution  is  $ \mathcal{O}(\ln(1/\epsilon))$,
 which matches the {\em optimal}  iteration complexity  for strongly convex optimization  in the deterministic regime.
\red{The oracle complexity of Algorithm \ref{alg-sto-gradient_track} for making $\left\| \begin{pmatrix}
\mathbb{E}[\| \bar{x}(k) -  x^*  \|] \\
\mathbb{E}[\| x(k) - (\mathbf{1} \otimes I_d)\bar{x}(k)\|]  \\
\mathbb{E}[ \| y(k) - (\mathbf{1} \otimes I_d)  \bar{y}(k)\|]
\end{pmatrix} \right\|\leq \epsilon$  is $\mathcal{O}(1/\epsilon^2)$    when   $q\in (\rho(J(\bm{\alpha})),1).  $
Recall that for the centralized SGD,
 the oracle complexity for making either the suboptimality gap  $\mathbb{E}[F(x)]-F(x^*)<\epsilon$
or the mean-squared error   $\mathbb{E}[\|x-x^*\|^2]<\epsilon$ is $O(1/\epsilon)$ (see e.g.,\cite{hazan2014beyond}),
which implies that the oracle complexity for obtaining $\mathbb{E}[\|x-x^*\|]<\epsilon$ is $O(1/\epsilon^2)$.
Thus, the number of sampled gradient required by   Algorithm \ref{alg-sto-gradient_track} with $N(k)=\lceil q^{-2k} \rceil,~q\in (\rho(J(\bm{\alpha})),1) $  to achieve a given solution accuracy   matches that of the centralized SGD.   \hfill $\Box$}
  %Therefore, the communication cost is saved without increasing the sample burden.

\end{remark}

\blue{Next, we investigate the communication complexity for obtaining an approximate solution.
We consider a special case with fixed graph and  impose the following condition.
 \begin{assumption}\label{ass-graph2}
 (i) $\mathcal{G}(k)\equiv \mathcal{G}$, where  $\mathcal{G}$ is strongly connected.

 (ii)     $A(k)\equiv A$, where the adjacency matrix $A$ associated with  $\mathcal{G}$ is doubly stochastic.

\end{assumption}
\begin{figure*}[hb]
\begin{minipage}[t]{0.3\linewidth}
 \includegraphics[width=2.4in]{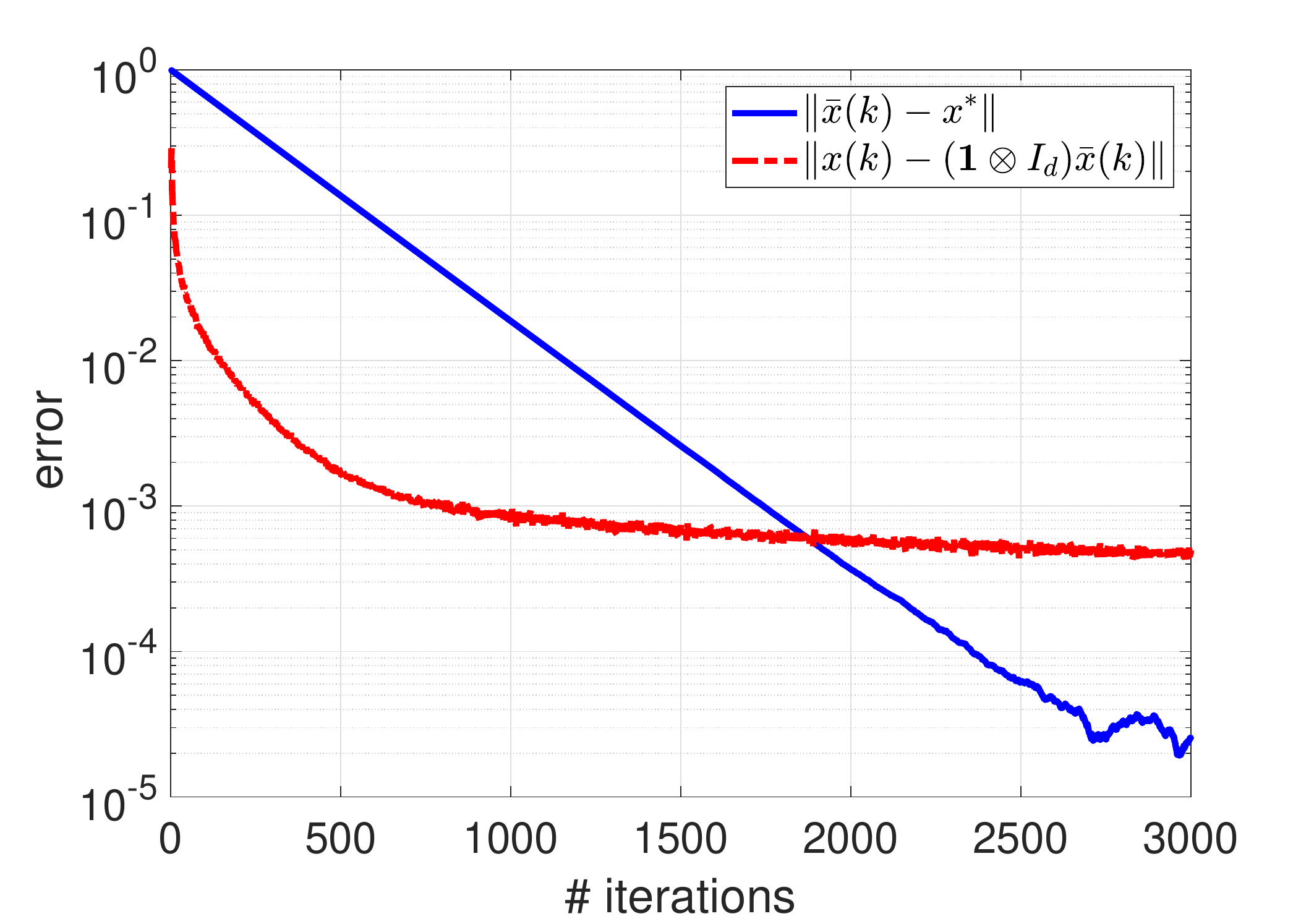}
 \caption{Convergence of Algorithm 1 for convex problems} \label{Gcon}
\end{minipage}
\begin{minipage}[t]{0.03\linewidth}
\end{minipage}
\begin{minipage}[t]{0.33\linewidth}
\centering \includegraphics[width=2.5in]{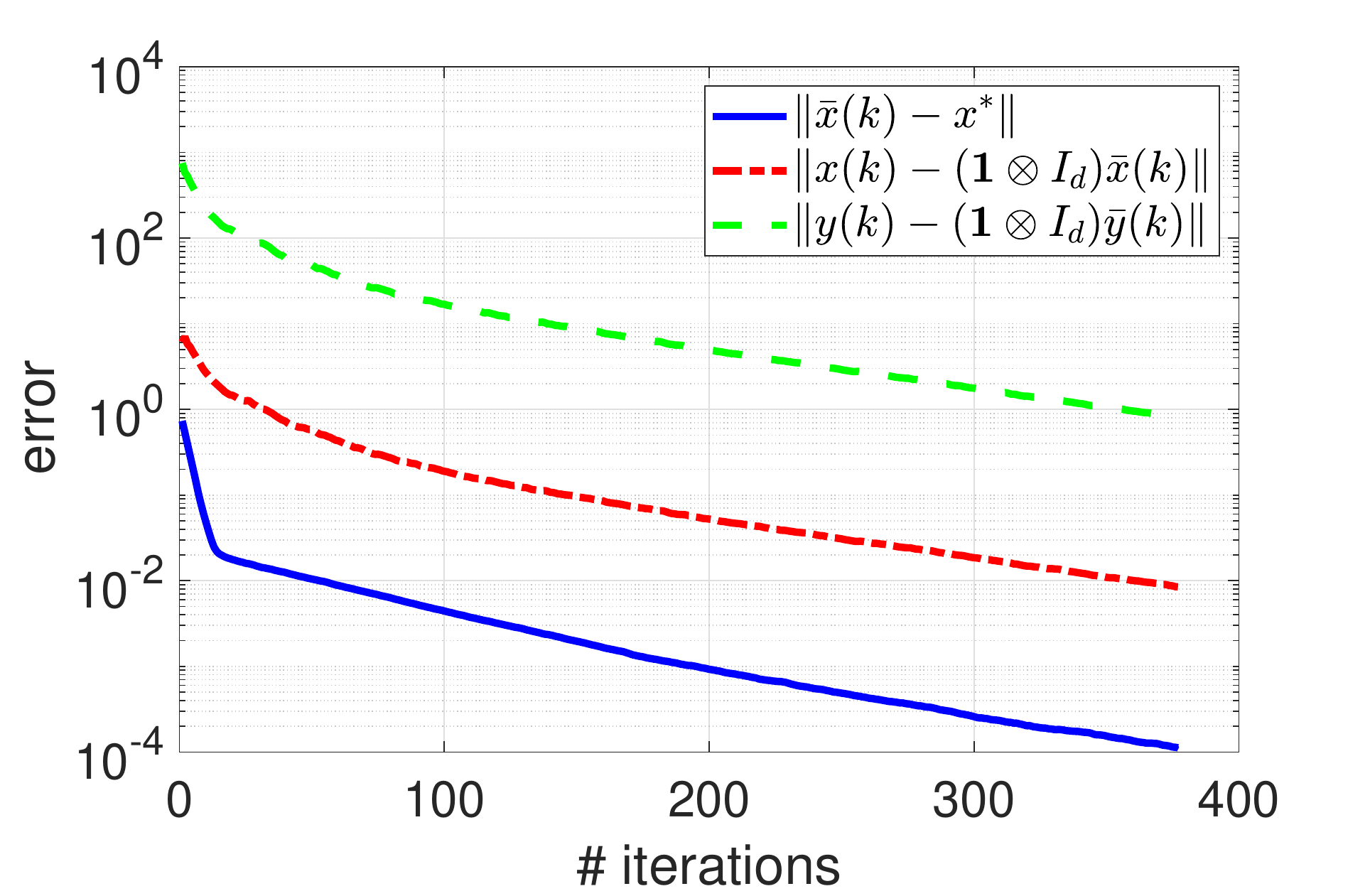}
\caption{Geometric  rate of Algorithm  1 for strongly convex problems }      \label{GOne}
\end{minipage}
\begin{minipage}[t]{0.03\linewidth}
\end{minipage}
\begin{minipage}[t]{0.3\linewidth}
\centering \includegraphics[width=2.3in]{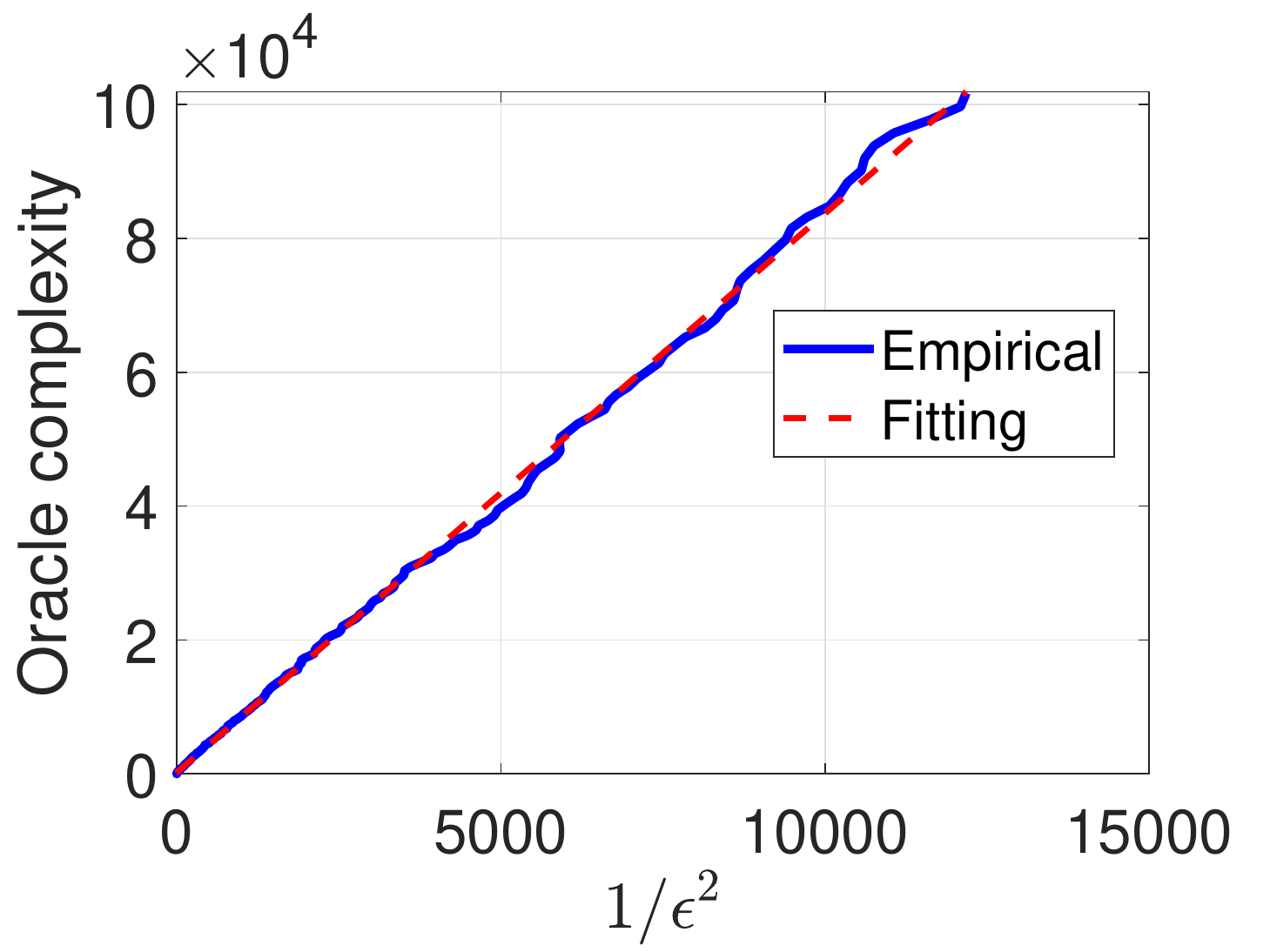}
  \caption{ Empirical oracle complexity and its  fitting of Algorithm  \ref{alg-sto-gradient_track}  for strongly convex problems }      \label{GThree}
\end{minipage}
\end{figure*}

   \begin{theorem}\label{thm3}  Let Assumptions   \ref{ass-convex}(ii), \ref{ass-noise}, \ref{ass-fun}, and \ref{ass-graph2} hold. Consider  Algorithm \ref{alg-sto-gradient_track} with   $N(k)=\lceil q^{-2k} \rceil,~q\in (0,1), $  where  the  step-size   $\alpha_i ,i\in \mathcal{V}$  satisfies \eqref{bd-stepsize}.
  Then  the number of communications required to obtain $\|z(k)\| \leq \epsilon$ is $\mathcal{O}(| \mathcal{E}|\ln(1/\epsilon)).$\end{theorem}
\begin{IEEEproof}
In each iteration $k$,  agent  $i$  requires  $2 | \mathcal{N}_i| $ rounds
of communication to obtain its neighbors' information $x_j(k) $ and $ y_j(k) .$
Thus, the number of  communication rounds required across the network at time $k$ is  $2| \mathcal{E}|.$
Since the number of iterations required to obtain $\|z(k)\| \leq \epsilon$ is $\mathcal{O}( \ln(1/\epsilon)),$ the  number of total communication rounds required is $ \mathcal{O}(|\mathcal{E}|\ln(1/\epsilon))$.
\end{IEEEproof}}

\blue{There might exist settings where a geometrically  increasing   batch-size  is impractical.  To this end, we
consider the use of polynomially  increasing batch-size that allows
for more gentle growth, and    proceed to investigate the  convergence rate as well as the complexity bounds.

\begin{theorem}\label{thm4}  Let Assumptions   \ref{ass-convex}(ii), \ref{ass-noise}, \ref{ass-fun}, and \ref{ass-graph2} hold. Consider  Algorithm \ref{alg-sto-gradient_track} with   $N(k)=\lceil (k+1)^{2\theta} \rceil,~\theta>0, $
  and the  step-size   $\alpha_i ,i\in \mathcal{V}$  satisfying  \eqref{bd-stepsize}.
  Then  $z(k)$   converges    to zero at a  polynomial   rate $\mathcal{O}(k^{-\theta})$. In addition, the number  of    samples and communications required to make  $\|z(k)\| \leq \epsilon$   is
     $\mathcal{O} \left( (1/\epsilon)^{2+1/\theta}  \right)$ and $ \mathcal{O}(|\mathcal{E}|  (1/\epsilon)^{1/\theta})$, respectively.  \end{theorem}
  \begin{IEEEproof}The proof is given in Appendix \ref{app:B5}.\end{IEEEproof}}

\red{\begin{remark}
 Though an increasing  batch-size implies a  higher sampling and computation burden than SGD  with a single iteration,
 the proposed scheme can significantly reduce the communication burden    compared  with  \cite{jakovetic2018convergence,yuan2018optimal,sayin2017stochastic,nedic2016stochastic,pu2020distributed}.
 Thus, Algorithm  \ref{alg-sto-gradient_track} is superior in   many practical networks  especially in  wireless networks, where the communication cost is usually much  higher than   gradient computations.
Therefore, the proposed scheme can remarkably save the communication cost by fully exploiting the local computation resources. Theorem \ref{thm3} and  Theorem \ref{thm4}      characterize  the trade-off between communication costs and sampling  rate, where a higher sampling rate leads to a smaller communication burden. \hfill $\Box$
 \end{remark}}

 \section{Numerical Simulations}

In this section, we examine the empirical   performance of  Algorithm \ref{alg-sto-gradient_track}  on the  distributed parameter estimation problems.
%and  nonlinear least-squares problems to demonstrate the behaviour of.
%Throughout this section, the empirical mean error is based on averaging across 50 trajectories.

\subsection{Distributed parameter estimation problem}
   Consider a network of  $n$ spatially distributed  sensors that aim to  estimate an  unknown  $d $-dimensional parameter $x^*$.
   Each sensor  $i$ collects a set of scalar    measurements   $ \{ d_{i,p} \}$ generated by  the  following linear regression model corrupted with noises,
 \begin{equation}
    d_{i,p}=u_{i,p}^Tx^* +\nu_{i,p},\nonumber
 \end{equation}
  where   $u_{i,p} \in  \mathbb{R}^{d}$ is  the regression vector  accessible to  agent $i$
  and  $\nu_{i,p}  \in \mathbb{R} $ is a zero-mean Gaussian noise.

Suppose that $\{u_{i,p}\}$ and $\{\nu_{i,p}\}$ are mutually independent  Gaussian sequences with  distributions $N(\mathbf{0}, R_{u,i})$ and    $N(0,\sigma_{i,\nu}^2)$, respectively.
Then the distributed parameter estimation problem can be modelled as a distributed stochastic  quadratic optimization problem,
 \begin{equation}\label{filter1}
  \min_{x\in \mathbb{R}^d}
  ~  {1\over n}\sum_{i=1}^n f_i(x) ,  {\rm~where~ }f_i(x)= \mathbb{E}   \big[\|  d_{i,p}-u_{i,p}^Tx\| ^2\big].
  \end{equation}
   Thus, $f_i(x)=(x-x^*)^T R_{u,i}(x-x^*)+\sigma_{i,\nu}^2 $ is convex and $\nabla f_i(x)= R_{u,i}(x-x^*).$
By using  the observed  regressor   $u_{i,p}$ and  the corresponding measurement $d_{i,p}$,
 the sampled gradient $ u_{i,p}u_{i,p}^T x-d_{i,p} u_{i,p}$   satisfies  Assumption \ref{ass-noise}.

\begin{figure*}[htbp]
\begin{minipage}[t]{0.3\linewidth}
\includegraphics[width=2.2in]{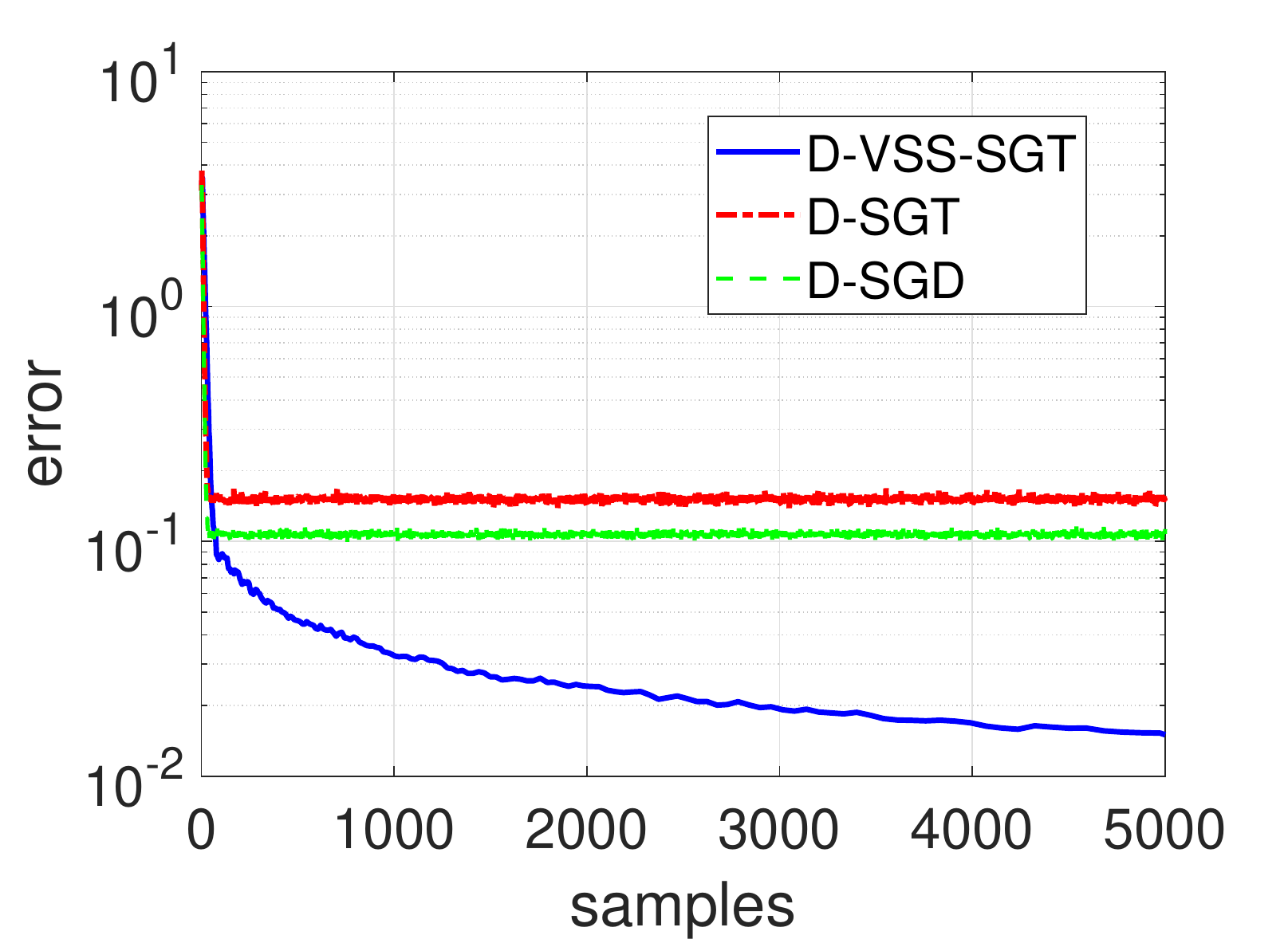}
 \caption{Comparison of Algorithm  \ref{alg-sto-gradient_track}  with  D-SGD and D-SGT under the same constant stepsize  }     \label{Giter}
\end{minipage}
\begin{minipage}[t]{0.03\linewidth}
\end{minipage}
\begin{minipage}[t]{0.32\linewidth}
 \includegraphics[width=2.2 in]{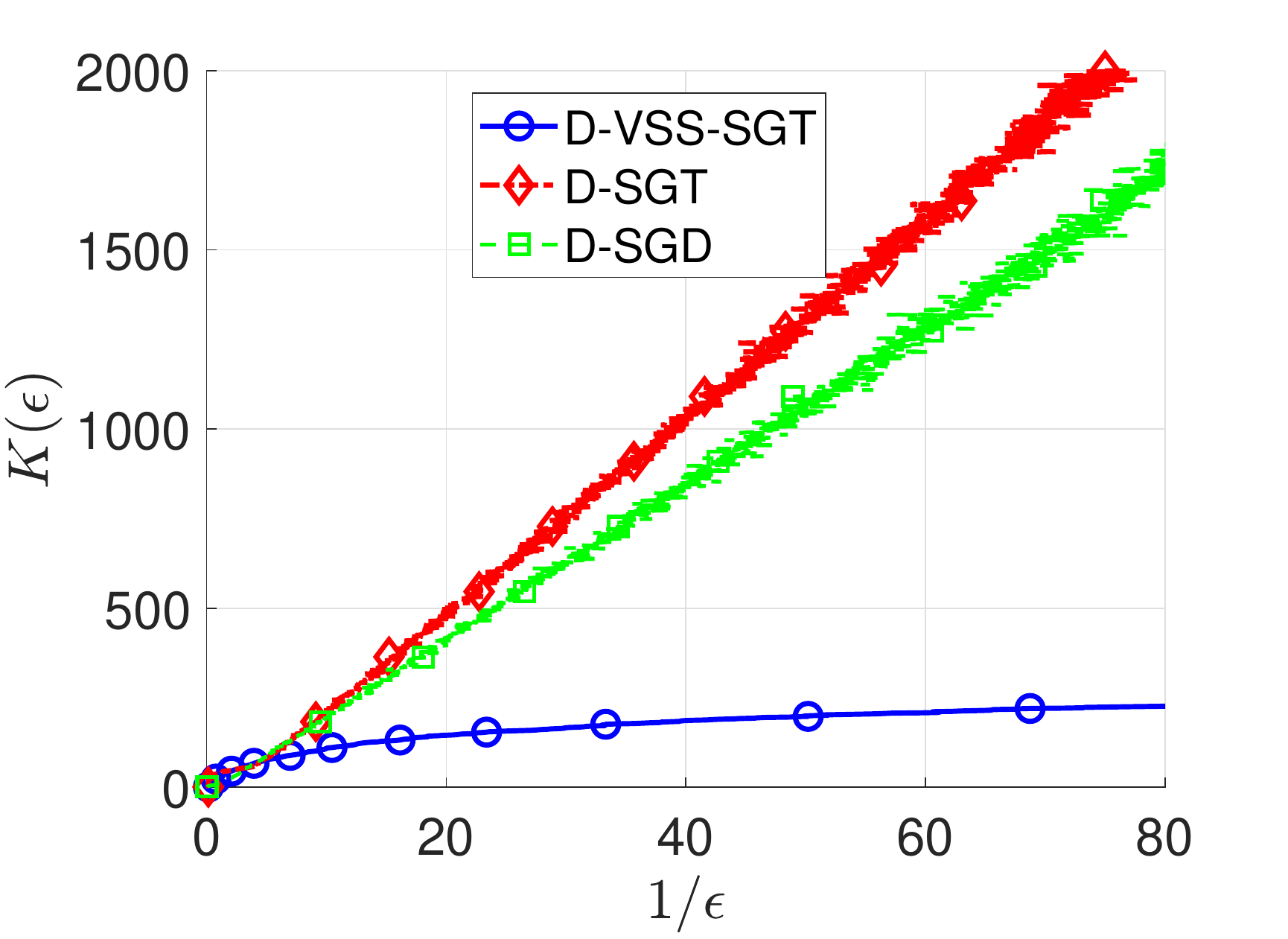}
 \caption{Iteration complexity  of  Algorithm  \ref{alg-sto-gradient_track}, D-SGD, and D-SGT}     \label{Gietr}
\end{minipage}
\begin{minipage}[t]{0.03\linewidth}
\end{minipage}
\begin{minipage}[t]{0.31\linewidth}
 \includegraphics[width=2.2 in]{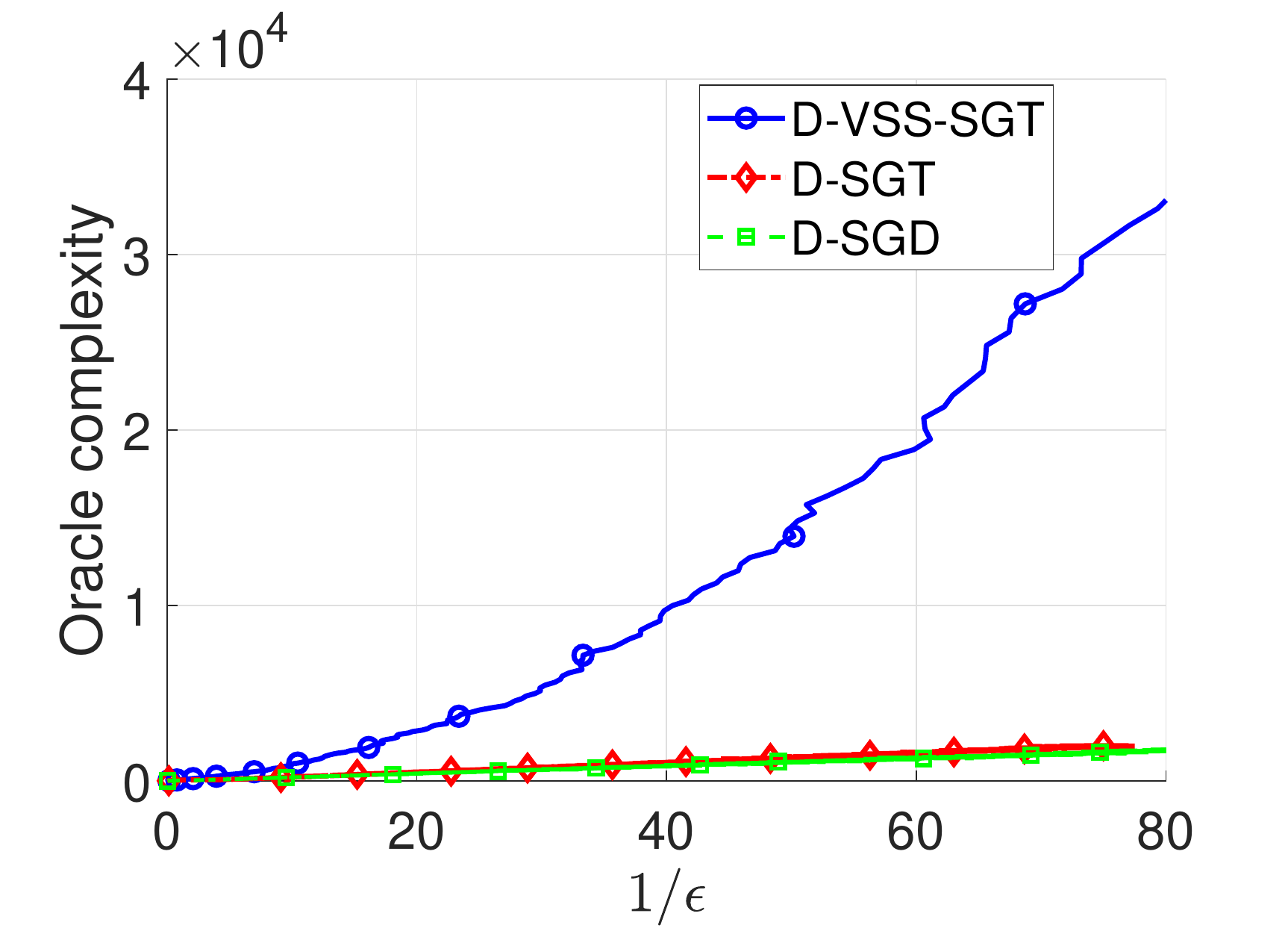}
 \caption{Oracle complexity of Algorithm  \ref{alg-sto-gradient_track}, D-SGD and D-SGT}     \label{Goracle}
\end{minipage}
\end{figure*}

\begin{figure*}[htbp]
\begin{minipage}[t]{0.3\linewidth}
 \includegraphics[width=2.2 in]{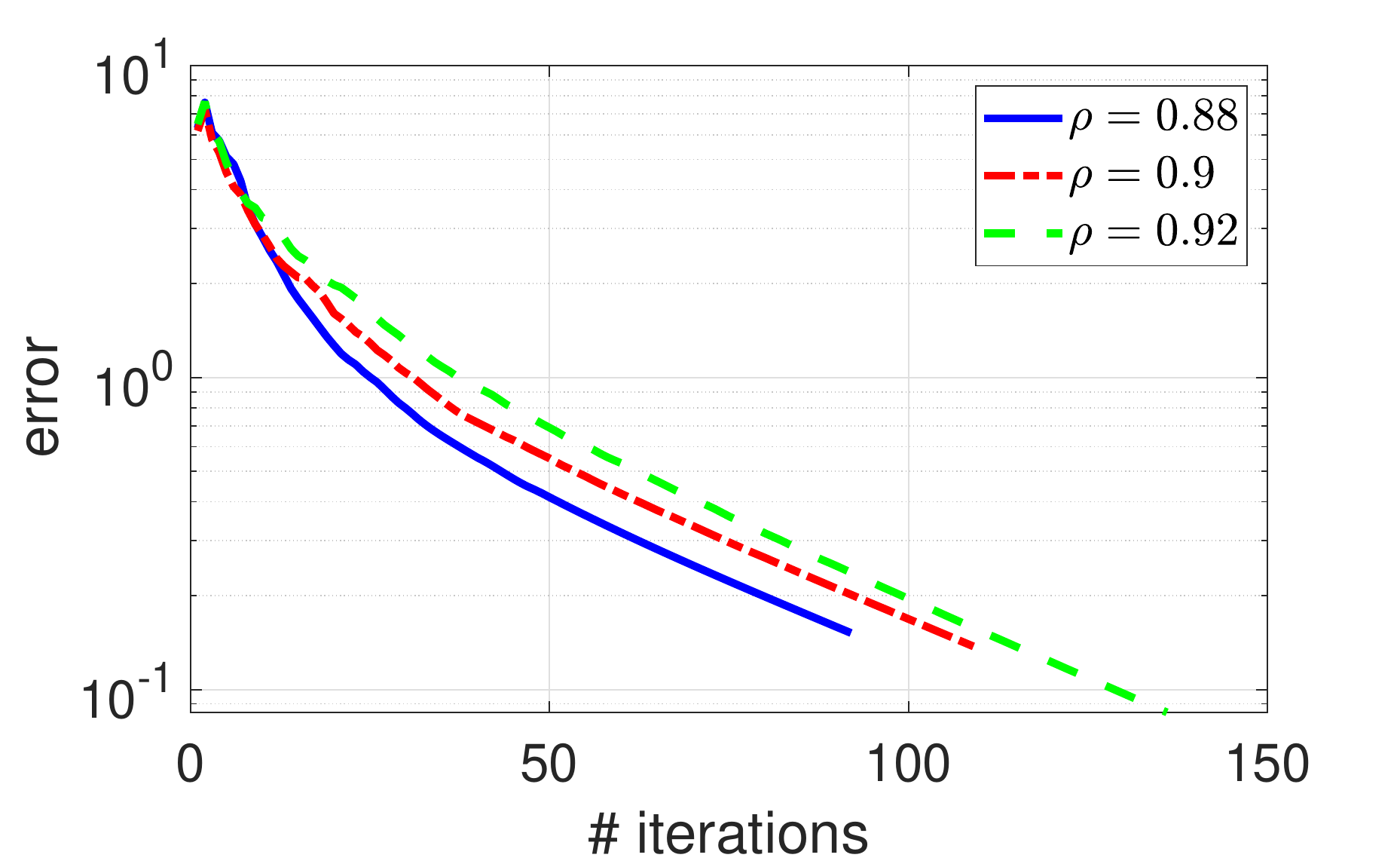}
 \caption{ Rate of $e(k)$  for  batch-sizes $N_k=\lceil \rho^{-k}\rceil $}     \label{G4}
\end{minipage}
\begin{minipage}[t]{0.03\linewidth}
\end{minipage}
\begin{minipage}[t]{0.33\linewidth}
\centering \includegraphics[width=2.1in]{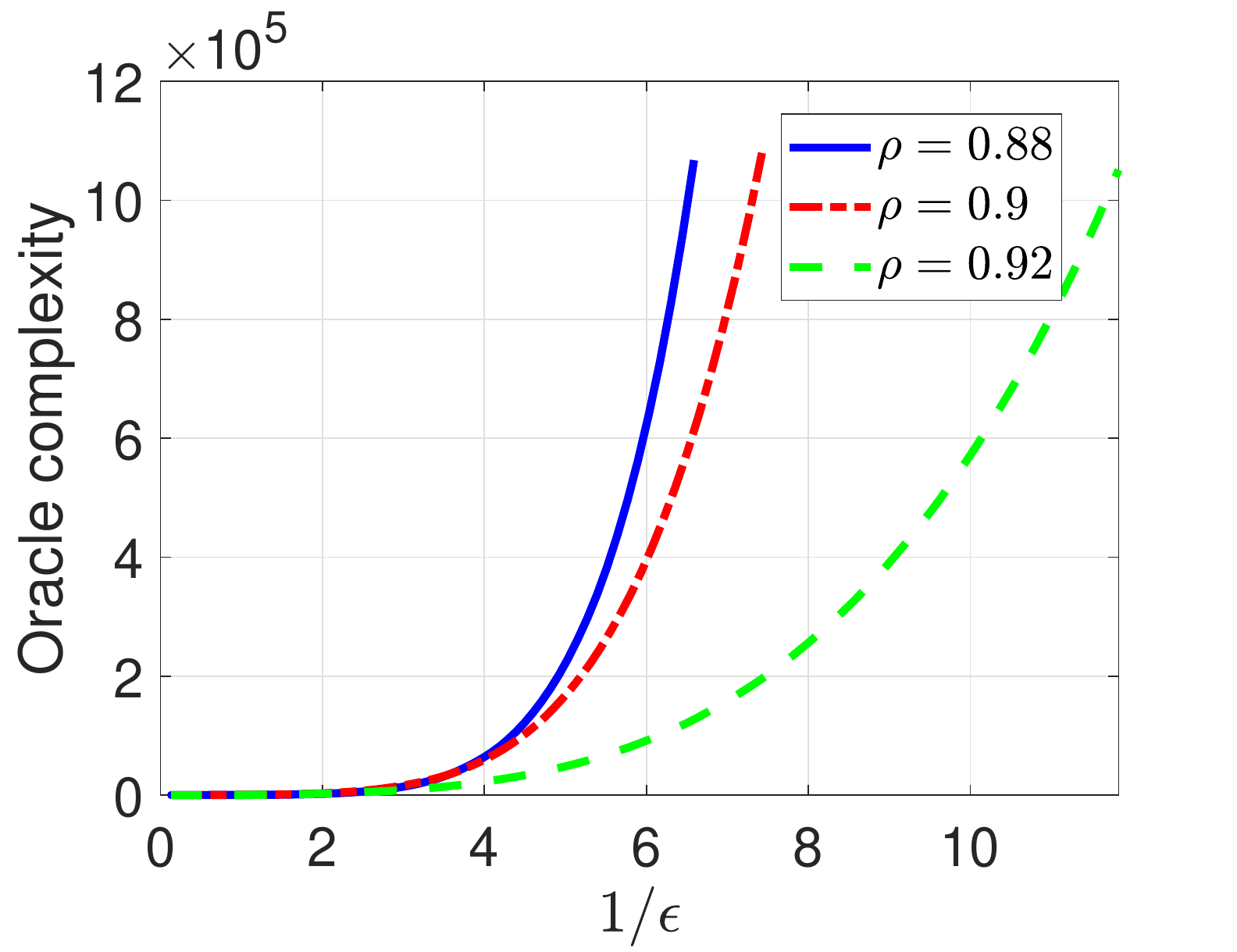}
 \caption{ Oracle complexity for   batch-sizes $N_k=\lceil \rho^{-k}\rceil $}     \label{G5}
\end{minipage}
\begin{minipage}[t]{0.03\linewidth}
\end{minipage}
\begin{minipage}[t]{0.3\linewidth}
 \includegraphics[width=2.2in]{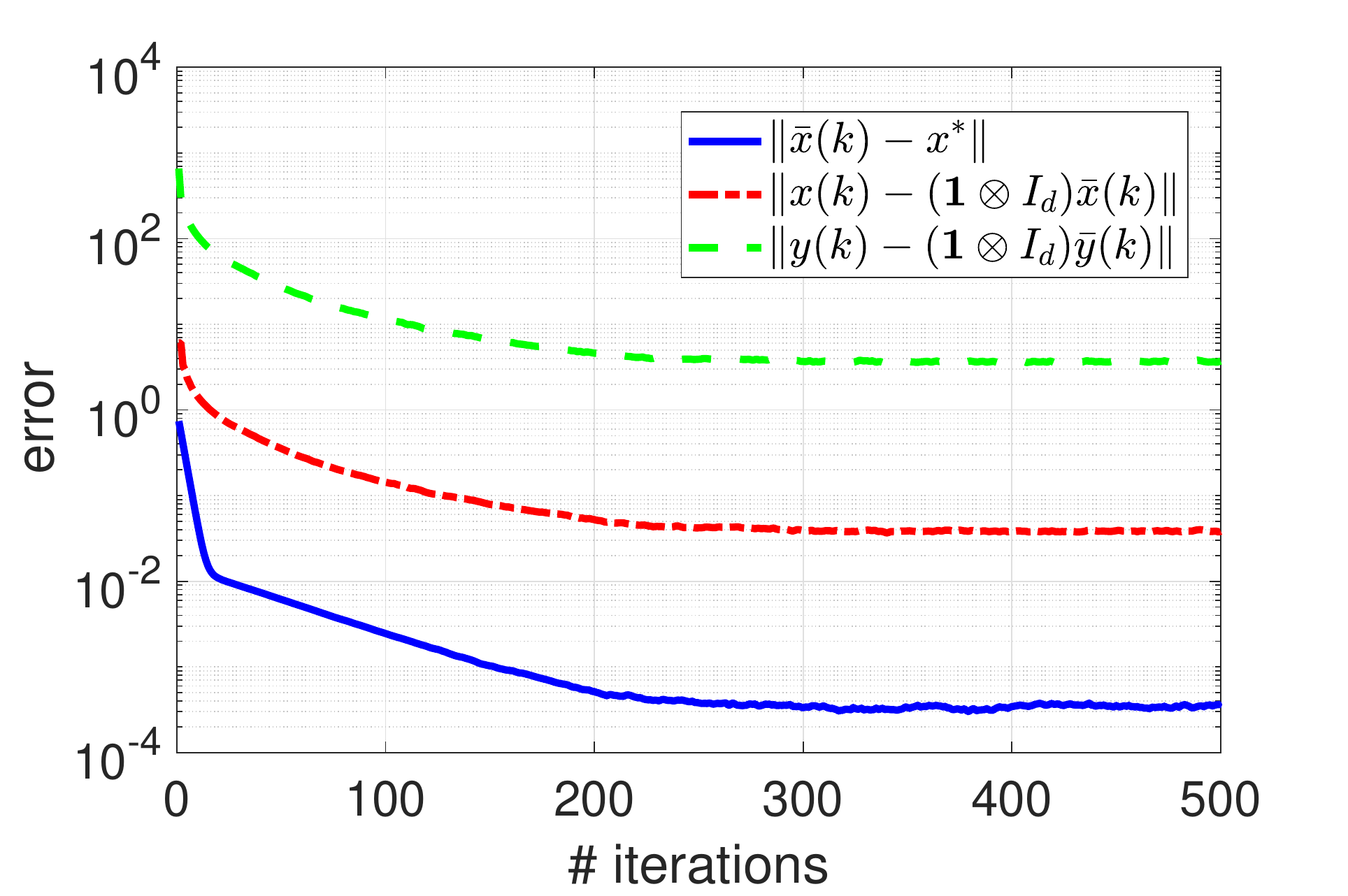}
 \caption{Algorithm performance with  constant   batch-size }     \label{G6}
\end{minipage}
\end{figure*}

%%%%%%%%%%%%%%%%%%%%%%%%%%%%%%%%%%%%%%%%%%%%%%%%%
\subsection{Numerical simulations}
%%%%%%%%%%%%%%%%%%%%%%%%%%%%%%%%%%%%%%%%%%%%%%%%%%

Set the vector dimension $d=10,$ the agent number $n=200,$ and the true parameter $x^*=\mathbf{1}/\sqrt{d} $.
We  randomly generate 10  undirected networks,  where any two  distinct agents are  linked with probability $0.1.$
The adjacency matrix is constructed based on the Metropolis rule.
A graph is uniformly sampled at each iteration such that Assumption \ref{ass-graph} is ensured.

{\bf Validation of Theorem  \ref{thm-as}.} Let each  covariance matrix $ R_{u,i}$  be  {\it positive semidefinite} with two eigenvalues equal to 0, that is,  each  $f_i(x)$ is {\it merely convex}. We run Algorithm \ref{alg-sto-gradient_track} with $\alpha=0.01$ and $N_k=\lceil k^{1.1}\rceil $,  and display the estimation errors of a sample path in
Fig. \ref{Gcon}, which shows that the generated iterates will asymptotically converge to the true parameter $x^*.$

%\begin{figure}[htbp]
%\begin{minipage}[t]{0.47\linewidth}
% \includegraphics[width=1.55in]{pic/rate_convex}
% \caption{Convergence of Alg. 1 for convex problems} \label{Gcon}
%\end{minipage}
%\begin{minipage}[t]{0.03\linewidth}
%\end{minipage}
%\begin{minipage}[t]{0.5\linewidth}
%\centering \includegraphics[width=1.55in]{pic/rate}
%\caption{Geometric  rate of Alg. 1 for strongly convex problems }      \label{GOne}
%\end{minipage}
%\end{figure}

%\begin{figure}[!htb]
%\centering
% \includegraphics[width=2.8in]{pic/rate_convex}
% \caption{Convergence of Alg. 1 for convex problems.} \label{Gcon}
% \vspace{-0.15in}
%\end{figure}

{\bf Validation of Theorems \ref{thm1} and \ref{thm2}.}
Let each  covariance matrix $ R_{u,i}$   be  positive definite. Then each  $f_i(x)$ is strongly convex and $x^*$ is the unique optimal solution to  \eqref{filter1}.  We run Algorithm \ref{alg-sto-gradient_track} with $\alpha=0.01$ and $N_k=\lceil 0.98^{-k}\rceil $, and
 examine the   empirical  rate  of convergence  and oracle complexity, where the empirical mean-squared error is based on averaging across 50 sample trajectories.
 The convergence rate    shown in Fig.  \ref{GOne},  demonstrating   that the iterates
$\{x_k\}$   generated by   Algorithm  \ref{alg-sto-gradient_track} converge to $x^*$ at a geometric  rate.  The oracle complexity   is shown in   Fig.  \ref{GThree}, where   x-axis is $1/\epsilon^2$
and   y-axis  denotes  the   number of    sampled gradients required to ensure
$$ e(k)\triangleq \mathbb{E}\left[\left\| \begin{pmatrix}  \bar{x}(k) -x^* \\ x(k) - (\mathbf{1}_n \otimes I_d)\bar{x}(k)  \end{pmatrix} \right \| \right]   <\epsilon.$$
 In Fig.  \ref{GThree}, the  blue solid curve represents the empirical data, while the  red    dashed curve denotes its linear fitting, which implies that the empirical
oracle complexity fits well with the established theoretical bound $\mathcal{O}(1/\epsilon^2)$.

%
%\begin{figure}[htbp]
%
%\begin{minipage}[t]{0.49\linewidth}
%\centering \includegraphics[width=1.5in]{pic/oracle}
%  \caption{ Empirical oracle complexity and its  fitting of Alg. \ref{alg-sto-gradient_track}  for strongly convex problems }      \label{GThree}
%\end{minipage}
%\begin{minipage}[t]{0.02\linewidth}
%\end{minipage}
%\begin{minipage}[t]{0.3\linewidth}
%\includegraphics[width=1.5in]{pic/compare}
% \caption{Comparison of Alg. \ref{alg-sto-gradient_track}  with  D-SGD and D-SGT under the same constant stepsize  }     \label{Giter}
%\end{minipage}
%\end{figure}

{\bf Comparison with \cite{ram2010distributed} and \cite{pu2020distributed}.}
We  compare  the performance of Algorithm  \ref{alg-sto-gradient_track}, abbreviated as D-VSS-SGT, with the  distributed stochastic gradient descent (D-SGD) \cite{ram2010distributed} and the distributed stochastic gradient tracking  (D-SGT) \cite{pu2020distributed} for strongly convex stochastic optimization.

Firstly, we  compare the algorithm performance of  the three methods  under fixed step-sizes.
We set    $\alpha=0.005$ in all  three  schemes, and   $N_k=\lceil 0.98^{-k}\rceil $  in Algorithm \ref{alg-sto-gradient_track}. The empirical error   $e(k)\triangleq \mathbb{E}\left[\left\| \begin{pmatrix}  \bar{x}(k) -x^* \\ x(k) - (\mathbf{1}_n \otimes I_d)\bar{x}(k)  \end{pmatrix} \right \| \right]  $  with respect to  the   number of sampled gradients  is given in Figure \ref{Giter}. It shows that the iterates of  D-SGD  and D-SGT  ceased at a neighborhood of the   true parameter $x^*$, while  the iterates generated by Algorithm  \ref{alg-sto-gradient_track}  will asymptotically converge   to  the true value  $x^*$.
It also shows that D-SGD  and D-SGT are more efficient in utilizing the samples than Algorithm  \ref{alg-sto-gradient_track} at  the first few samples, but with the increasing of gradient samples, Algorithm  \ref{alg-sto-gradient_track}  is
\red{superior} than D-SGD  and D-SGT.

We further compare the iteration and oracle complexity of the three methods, where Algorithm  \ref{alg-sto-gradient_track} uses a constant stepsize,  while D-SGD   and D-SGT   use decreasing stepsizes. The empirical number of iterations
 and sampled gradients required to obtain a solution with
 the same accuracy are demonstrated  in Fig.  \ref{Gietr} and Fig. \ref{Goracle}.
 We see from Fig.  \ref{Gietr}    that Algorithm  \ref{alg-sto-gradient_track}  can significantly reduce the iteration numbers, so do the   communication rounds (costs). Meanwhile,  Fig. \ref{Goracle} shows that
 Algorithm  \ref{alg-sto-gradient_track} requires more sampled gradients than D-SGD and D-SGT.
 In fact, in certain applications such as wireless networks, high communication overhead may render a distributed scheme impractical.
  As such, the variance-reduced method  proposed in this work is suitable for \red{network problems} when the communication costs are more expensive than sampling and local computations.

%\begin{figure}[htbp]
%\begin{minipage}[t]{0.3\linewidth}
% \includegraphics[width=2.2 in]{pic/iter_three}
% \caption{Iteration complexity  of  Alg. \ref{alg-sto-gradient_track}, D-SGD, and D-SGT}     \label{Gietr}
%\end{minipage}
%\begin{minipage}[t]{0.02\linewidth}
%\end{minipage}
%\begin{minipage}[t]{0.3\linewidth}
% \includegraphics[width=2.2 in]{pic/oracle_three}
% \caption{Oracle complexity of Alg. \ref{alg-sto-gradient_track}, D-SGD and D-SGT}     \label{Goracle}
%\end{minipage}
%\end{figure}

%
%\begin{figure}[!htb]
%\centering
% \includegraphics[width=2.8 in]{pic/iter_three}
% \caption{Iteration complexity of  Alg. \ref{alg-sto-gradient_track}, D-SGD, and D-SGT}     \label{Gietr}
%\vspace{-0.15in} \end{figure}
%
%
%\begin{figure}[!htb]
%\centering
% \includegraphics[width=2.8 in]{pic/oracle_three}
% \caption{Oracle complexity of Alg. \ref{alg-sto-gradient_track}, D-SGD and D-SGT}     \label{Goracle}
%\vspace{-0.15in} \end{figure}

{\bf Influence of the   batch-sizes.}  We run Algorithm \ref{alg-sto-gradient_track}  with $\alpha=0.01$ and different  geometric batch-sizes $N_k=\lceil \rho^{-k}\rceil $. We set $\rho=0.88,0.9,0.92$,
  and display the empirical rate and oracle complexity in  Fig. \ref{G4}
 and Fig. \ref{G5}, respectively. We conclude from the figures  that a faster increasing batch-size   leads to a better convergence rate (namely,  less rounds of communications)  while at the  cost of more sampled data and heavier   computations. Hence, the parameter $\rho$ should be properly selected to balance
 communication costs, sampling costs and computation costs in practice.

%
%\begin{figure}[htbp]
%\begin{minipage}[t]{0.49\linewidth}
% \includegraphics[width=1.8 in]{pic/rate_sample}
% \caption{ Rate of $e(k)$  for  batch-sizes $N_k=\lceil \rho^{-k}\rceil $}     \label{G4}
%\end{minipage}
%\begin{minipage}[t]{0.04\linewidth}
%\end{minipage}
%\begin{minipage}[t]{0.48\linewidth}
%  \includegraphics[width=1.45in]{pic/oracle_sample}
% \caption{ Oracle complexity for   batch-sizes $N_k=\lceil \rho^{-k}\rceil $}     \label{G5}
%\end{minipage}
%\end{figure}

\blue{{\bf Performance of Algorithm \ref{alg-sto-gradient_track} with  constant batch-size.}
 Finally, we run the algorithm   with $\alpha=0.01$ and a  constant batch-size  $N_k\equiv 20  $,
 and display the empirical convergence rate  in  Fig. \ref{G6}, which  clearly shows  that  the
  estimates ceased at a neighborhood of the optimal solution.}

%\begin{figure}[!htb]
%\centering
% \includegraphics[width=2.5in]{pic/iter_constant}
% \caption{Algorithm performance with  constant   batch-size }     \label{G6}
%\vspace{-0.15in} \end{figure}
\section{Conclusions}

We  proposed a distributed stochastic gradient tracking algorithm with variable sample-sizes for stochastic optimization over  random networks.  We proved that with  a  suitably selected  constant stepsize  and a properly increasing gradient sample-size,
the iterates converge almost surely to the optimal solution for convex problmes.
 For strongly convex problems, we further obtained  the geometric convergence rate  with geometrically increasing batch-sizes and  established the complexity bounds for obtaining an $\epsilon$-optimal solution.
Both the iteration complexity and the oracle complexity  are comparable with  the centralized stochastic gradient descent algorithm. \red{It   might be of interests  to embed   the push-pull method  for resolving distributed stochastic optimization with general digraphs.}
The extension of the current algorithm to non-convex/non-smooth distributed stochastic optimization
is a promising future research direction.

\appendices
\numberwithin{equation}{section}

\section{Proofs of  Section \ref{sec:as}}\label{app:A}
\subsection{  Proof of Lemma \ref{lem-cond}.} \label{app:A1}
  {\em We  first give a bound on $ \mathbb{E}[\| \tilde{x} (k+1)  \| | \mathcal{F}(k)].$}
Note by Assumption \ref{ass-graph}(i) and $D_{\bot}=I_n-{\mathbf{1}_n\mathbf{1}_n^T / n} $ that $D_{\bot}A(k)=(A(k)- \mathbf{1}_n\mathbf{1}_n^T / n ) D_{\bot}$.   Then by multiplying both sides of  \eqref{push2} with $D_{\bot}\otimes I_d$  from the left, using the definition \eqref{def-tidlex}, and $\bar{\alpha}={\sum_{i=1}^n \alpha_i \over n}$, we obtain that
\begin{align*}
 &\tilde{x} (k+1) =  (A(k)-\mathbf{1}_n\mathbf{1}_n^T / n) \otimes I_d  \tilde{x}(k)   - (\bm{\alpha}  \otimes I_d)  y(k)
\\  & \qquad +    {\mathbf{1}_n \mathbf{1}_n^T \bm{\alpha} \otimes I_d \over n} y (k)\\
& \overset{\eqref{def-tidley}}{ =}  (A(k)-\mathbf{1}_n\mathbf{1}_n^T / n) \otimes I_d  \tilde{x}(k)   - (\bm{\alpha}  \otimes I_d) \tilde{ y}(k)
\\ &    -(\bm{\alpha} \mathbf{1}_n \otimes I_d)   \bar{y}(k)+  \bar{\alpha}  {\mathbf{1}_n   \otimes I_d  }\bar{y}(k)+ {\mathbf{1}_n  \mathbf{1}_n^T \bm{\alpha} \otimes I_d \over n} \tilde{y}(k) .
\end{align*}
Since  $\bar{y}(k), \tilde{y}(k)$ are adapted to  $\mathcal{F}(k)$,  by using the triangle  inequality,  we obtain that
\begin{equation*}
\begin{split}
 &\mathbb{E}[\| \tilde{x} (k+1)  \| | \mathcal{F}(k)] \leq  \mathbb{E}\big[\| (A(k)-\mathbf{1}_n\mathbf{1}_n^T / n) \otimes I_d  \tilde{x}(k) \| | \mathcal{F}(k) \big]
\\ &   +  \| \bar{\alpha}  \mathbf{1}_n   - \bm{\alpha} \mathbf{1}_n   \|  \| \bar{y}(k)\| + \| (I_n- \mathbf{1}_n  \mathbf{1}_n^T /n) \otimes I_d  \bm{\alpha}   \| \| \tilde{y}(k)\|.
\end{split}
\end{equation*}
This  combined with   \eqref{bd-sr}, \eqref{def-alpha},   and  $ \| (I_n- \mathbf{1}_n  \mathbf{1}_n^T /n) \otimes I_d  \bm{\alpha}   \| \leq \| (I_n- \mathbf{1}_n  \mathbf{1}_n^T /n) \otimes I_d     \| \| \bm{\alpha}   \| \leq \|   \bm{\alpha}   \| =\max_i \alpha_i$  proves \eqref{recur-tildex}.
 \vskip 2mm
{\em Next, we give a bound on $ \mathbb{E}[  \| \tilde{y}(k+1)  \| | \mathcal{F}(k)].$}
From \eqref{push2} and $A(k)\mathbf{1}_n =\mathbf{1}_n $ it follows that
\begin{align}
&  \| x (k+1)-x(k) \|    \notag
%\\&  = \|(A(k)-I_n)  \otimes I_d   x(k)- (\bm{\alpha}  \otimes I_d)  y(k) \| \notag
 \\& =\| (A(k)-I_n)  \otimes I_d  \tilde{x}(k)  \notag
\\& \quad - (\bm{\alpha}  \otimes I_d) \tilde{ y}(k)-(\bm{\alpha} \mathbf{1}_n \otimes I_d)   \bar{y}(k) \|  \notag
\\& \leq     \| \tilde{x}(k) \|  +\alpha_{\max}  \|   \tilde{y}(k)  \|   +c_2 \|    \bar{y}(k)  \| ,  \label{recur-xdiff}
 \end{align}
where in the last inequality  we use the triangle inequality, $\|  A(k)-I_n \|  \leq 1,$   $\|\bm{\alpha}  \|  = \alpha_{\max}
 $,  and $\|  \bm{\alpha} \mathbf{1}_n  \|=c_2 . $
Note by  \eqref{def-fw} and Assumption~\ref{ass-convex}(ii) that
\begin{align}
&\|  \red{\nabla f (k+1)-\nabla f(k) } \| \notag
\\& =\sqrt{\sum_{i=1}^n \|\nabla f_i(x_i(k+1))-\nabla f_i(x_i(k)) \|^2}\label{recur-gdiff}
\\&  \leq  \sqrt{\sum_{i=1}^nL^2  \| x_i(k+1)- x_i(k) \|^2} = L  \| x(k+1)- x(k) \| . \notag
\end{align}

By multiplying both sides of  \eqref{pull2} with $D_{\bot}\otimes I_d$
  from the left, and  using the definition \eqref{def-tidley},  we obtain that
  \begin{align*}
 &\tilde{y} (k+1) =  (A(k)-\mathbf{1}_n\mathbf{1}_n^T / n) \otimes I_d  \tilde{y}(k)
 \\&+ D_{\bot}\otimes I_d ( \red{\nabla f (k+1)-\nabla f(k) }+w(k+1)-w(k )).
\end{align*}
Then by the triangle  inequality, using \eqref{recur-xdiff}, \eqref{recur-gdiff}, and   $\| D_{\bot}\|\leq 1$,
we obtain that
 \begin{equation}\label{inequl-tdx}
\begin{split}
 &\| \tilde{y} (k+1) \leq \|(A(k)-\mathbf{1}_n\mathbf{1}_n^T / n) \otimes I_d  \tilde{y}(k)\|
 \\&+L  \| x(k+1)- x(k) \|+ \|w(k+1)-w(k ))\|
 \\& \leq  \|(A(k)-\mathbf{1}_n\mathbf{1}_n^T / n) \otimes I_d  \tilde{y}(k)\|+
  \alpha_{\max} L   \| \tilde{y}(k)\|
 \\&+L  \| \tilde{x}(k) \|     +c_2 L \|    \bar{y}(k)  \|  + \|w(k+1)-w(k ))\|.
 \end{split}
\end{equation}

 Since  $ \tilde{y}(k)$ is adapted to  $\mathcal{F}(k)$, and  $A(k)$ is
 independent of $\mathcal{F}(k) $, similarly to  \eqref{bd-sr}, we can obtain that
 \begin{equation*}
\begin{split}
&  \mathbb{E}\big[\| (A(k)-\mathbf{1}_n\mathbf{1}_n^T / n) \otimes I_d\tilde{y}(k)\| | \mathcal{F}(k) \big]\leq \rho_1 \| \tilde{y}(k)\|  .
\end{split}
\end{equation*}
Then by taking conditional expectation of \eqref{inequl-tdx} on  $\mathcal{F}(k),$  and recalling that  $ \tilde{y}(k),   \bar{y}(k),  \tilde{x}(k)$ are adapted to $\mathcal{F}(k) $,  we prove \eqref{recur-yerr2}.
\hfill $\Box$

\subsection{ Proof of Lemma \ref{lem-exp}.}\label{app:A2}

Since  for each  $i\in \mathcal{V}$, the   samples $\{ \xi^p_i(k)\}_{p=1}^{N(k)}$ are    independent, by Assumption \ref{ass-noise} and the definition \eqref{def-w}, we have
\begin{equation*}
\begin{split}
& \mathbb{E}[\| w_i(k) \|^2 |x_i(k) ]
\\& ={1\over N(k)^2} \sum_{p=1} ^{N(k) } \mathbb{E}   [  \|   \nabla  h_i(x_i(k),\xi^p_i(k))- \nabla  f_i(x_i(k))    \|^2 |x_i(k)   ] .
\\& \leq {\nu^2 \over N(k)}, \quad \forall k\geq 0.
\end{split}
\end{equation*}
This implies that
\begin{align}\label{bd-noisevr}
\mathbb{E}[\|w(k)\|^2]=\sum_{i=1}^n \mathbb{E}[\|w_i(k)\|^2]\leq {n\nu^2 \over N(k)}.
\end{align}
Then by  the Jensen's inequliaty,  we obtain that for any $k\geq 0$,
\begin{equation}\label{bd-noise}
\begin{split}  \mathbb{E}[\| w (k) \|] &\leq\sqrt{\mathbb{E}[\|w(k)\|^2] }  \leq
{ \sqrt{n} \nu  \over \sqrt{N(k)}}.
\end{split}
\end{equation}

By \eqref{bd-noise} and the triangle inequality, there holds $ \mathbb{E}[ \| w(k+1)-w(k ) \| ] \leq  p_k.$
By taking the unconditional expectations on both sides of \eqref{recur-tildex} and \eqref{recur-yerr2}, we obtain that   for any $k\geq 0 $,
   \begin{align}\label{bd-xy}
     & \mathbb{E} \left(
     \begin{array}{c}  \| \tilde{x} (k+1)  \|   \\
      \| \tilde{y} (k+1)  \|    \\
     \end{array}
   \right) \leq \underbrace{\left(
                  \begin{array}{cc}
                    \rho_1 & \alpha_{\max}\\
                    L &  \rho_1+\alpha_{\max}L \\
                  \end{array}
                \right)}_{\triangleq \Lambda} \mathbb{E} \left(
     \begin{array}{c}  \| \tilde{x} (k )  \|   \\
      \| \tilde{y} (k )  \|    \\
     \end{array}
   \right) \notag\\ &+\left(
     \begin{array}{c} c_1 \\
    c_2L   \\
     \end{array}
   \right)\mathbb{E}[\| \bar{y}(k)\|] +   \left(
                                          \begin{array}{c}
                                            0 \\
                                       p_k \\
                                          \end{array}
                                        \right).
\end{align}
Clearly,  the spectral radius  of the matrix $\Lambda$  is $\rho_2 ={2\rho_1+\alpha_{\max}L +\sqrt{\alpha_{\max}^2L^2+4\alpha_{\max}L} \over 2}.$
  Note by  $\alpha_i< { (1-\rho_1)^2\over (2-\rho_1)L}$  that $\alpha_{\max}L<{ (1-\rho_1)^2\over (2-\rho_1)} ,$ which implies that $\rho_2<1.$

By \eqref{bd-xy} and the definition of $e(k),$  we obtain   that
 \begin{align*}
  &e(k+1)  \leq  \rho_2 e(k)+ \sqrt{ c_1^2+   c_2^2L^2}    \mathbb{E}[\| \bar{y}(k)\|]   +   p_k.
\end{align*}
Therefore, we recursively obtain that
\begin{align*}
e(k)& \leq \rho_2^ke(0)+   \sum_{t=0}^{k-1} \rho_2^t   p_{k-1-t}
    \\&+\sqrt{ c_1^2+   c_2^2L^2}   \sum_{t=0}^{k-1} \rho_2^t \mathbb{E}[\| \bar{y}(k-1-t)\|] .
\end{align*}
Taking the square on both sides of the above equation,  using $(a+b+c)^2\leq 3(a^2+b^2+c^2) $
and  the  Cauchy-Schwarz inequality, we have
\begin{align*}
e(k)^2& \leq 3 \rho_2^{2k} e(0)^2+ 3  \left (\sum_{t=0}^{k-1} \rho_2^t p_{k-1-t}\right)^2
    \\&+3( c_1^2+   c_2^2L^2) \left(  \sum_{t=0}^{k-1} \rho_2^t \mathbb{E}[\| \bar{y}(k-1-t)\|] \right)^2
    \\& \leq  3 \rho_2^{2k} e(0)^2+ 3  \sum_{t=0}^{k-1} (\rho_2^{t/2})^2\sum_{t=0}^{k-1} \left(\rho_2^{t/2}  p_{k-1-t} \right)^2
    \\&+3( c_1^2+   c_2^2L^2)  \sum_{t=0}^{k-1}  (\rho_2^{t/2})^2 \sum_{t=0}^{k-1}\left(\rho_2^{t/2} \mathbb{E}[\| \bar{y}(k-1-t)\|] \right)^2
     \\& \leq  3 \rho_2^{2k} e(0)^2+ { 3   \over 1-\rho_2}\sum_{t=0}^{k-1} \rho_2^{t }  p_{k-1-t}  ^2
    \\&+{3( c_1^2+   c_2^2L^2)\over 1-\rho_2}  \sum_{t=0}^{k-1}\rho_2^{t} \mathbb{E}[\| \bar{y}(k-1-t)\|] ^2.
\end{align*}
where the last inequality holds  by
$ \sum_{t=0}^{k-1} \rho_2^{t } \leq {1\over 1-\rho_2}$.  Note that
\begin{align*}
\sum_{k=1}^K \sum_{t=0}^{k-1} \rho_2^{t }  p_{k-1-t}^2=\sum_{s=0}^{K-1} p_s^2  \sum_{k=0}^{K-1-i} \rho_2^k
\leq {1\over 1-\rho_2} \sum_{s=0}^{K } p_s^2 .
\end{align*}
By summing the relation over $k$ from $1$ to $K$ and adding $e(0)^2$ to both sides, we obtain that
\begin{align*}
&\sum_{k=0}^K e(k)^2 \leq 3 \sum_{k=0}^K\rho_2^{2k} e(0)^2+   { 3    \over 1-\rho_2}\sum_{k=1}^K\sum_{t=0}^{k-1} \rho_2^{t }  p_{k-1-t}  ^2   \\&+{3( c_1^2+   c_2^2L^2)\over 1-\rho_2}  \sum_{k=1}^K \sum_{t=0}^{k-1}\rho_2^{t} \mathbb{E}[\| \bar{y}(k-1-t)\|] ^2
    \\& \leq {3 e(0)^2 \over 1-\rho_2^2}+ { 3  \over (1-\rho_2)^2}  \sum_{s=0}^{K } p_s^2+  {3( c_1^2+   c_2^2L^2)\over (1-\rho_2)^2} \sum_{s=0}^{K }\mathbb{E}[\| \bar{y}(s)\|] ^2.
\end{align*}
This combined  with $\sqrt{a^2+b^2+c^2}\leq a+b+c$ proves the lemma.
\hfill $\Box$

\subsection{Proof of Theorem \ref{thm-as}} \label{app:A3}
By  multiplying both sides of      \eqref{push2} with ${ (\mathbf{1}_n^T  \otimes I_d ) \over n} $ from the left
and using Assumption \ref{ass-graph}(i),  we   obtain  that
\begin{equation}\label{recur-barx0}
\begin{split}  & \bar{x}(k+1)-\bar{x}(k)= -{\mathbf{1}_n^T \bm{\alpha} \otimes I_d \over n} y (k) ,\quad \forall k\geq 0    \\
& \overset{ \eqref{def-tidley}}{=}  -\bar{\alpha}\bar{y}(k)- {(\mathbf{1}_n^T \bm{\alpha} - \bar{\alpha} \mathbf{1}_n^T)\otimes I_d \over n} \tilde{y} (k) .
\end{split}
\end{equation}
By the triangle inequality and \eqref{def-alpha}, we have
\begin{align}\label{bd-barx}  \| \bar{x}(k+1)-\bar{x}(k)\| & \leq \bar{\alpha} \| \bar{y}(k)\|+  {c_1  \over  n}    \|\tilde{y} (k)\| .
\end{align}
 Also, by using  \eqref{pull} and Assumption \ref{ass-graph}(i), we obtain that
\begin{equation*}
\begin{split}
\bar{y}(k+1) & =\bar{y}(k) +{1\over n} \sum_{i=1}^n \tilde{g}_i(x_i(k+1))
-  {1\over n} \sum_{i=1}^n \tilde{g}_i(x_i(k))  .
\end{split}
\end{equation*}
Then by recalling that   $y_i(0)= \tilde{g}_i(x_i(0))$,   one can recursively show that $
 \bar{y}(k) = {1\over n} \sum_{i=1}^n \tilde{g}_i(x_i(k)) $ for any $ k\geq 0.$
From   \eqref{def-w} it follows  that
 \begin{align}\label{recur-bary1}
 \bar{y}(k) = {1\over n} \sum_{i=1}^n \nabla f_i(x_i(k)) +{1\over n} \sum_{i=1}^n w_i(k) , \quad \forall k\geq 0.
\end{align}

 Denote
 \begin{align}\label{def-gk}   g(k)\triangleq \nabla  F(\bar{x}(k ))={1\over n} \sum_{i=1}^n  f_i (\bar{x}(k )).\end{align}
By Assumption \ref{ass-convex} and  the definition $F(x )= {1\over n } \sum_{i=1}^n f_i(x)$, we conclude that
 $F(x)$ is convex and its gradient function  is $L$-Lipschitz continuous. Therefore,
 \begin{equation*}
\begin{split}
 &F(\bar{x}(k+1))\leq  F(\bar{x}(k ))  +(\bar{x}(k+1)-\bar{x}(k))^T g(k)
  \\& +{L\over 2} \| \bar{x}(k+1)-\bar{x}(k) \|^2
 \\& \overset{\eqref{recur-bary1} }{ =} F(\bar{x}(k )) +{L\over 2} \| \bar{x}(k+1)-\bar{x}(k) \|^2
  \\&  + (\bar{x}(k+1)-\bar{x}(k))^T  \left( g(k) -{1\over n} \sum_{i=1}^n \nabla f_i(x_i(k)) \right)
  \\&+ (\bar{x}(k+1)-\bar{x}(k))^T \left(\bar{y}(k)-{1\over n} \sum_{i=1}^n w_i(k) \right) .
\end{split}
\end{equation*}

Define $v(k)\triangleq F(\bar{x}(k )) -F^*$.  Note by \eqref{recur-barx0} and the definition of $\mathcal{F}(k)$ that $\bar{x}(k+1)-\bar{x}(k)$ is adapted to $\mathcal{F}(k)$. By subtracting  $F^*$ from the above equation and taking the conditional expectation on $\mathcal{F}(k)$, we obtain that
 \begin{equation}\label{bd-Fx}
\begin{split}
 & \mathbb{E}[v(k+1)|\mathcal{F}(k)]\leq v(k)     + \underbrace{ {L\over 2} \| \bar{x}(k+1)-\bar{x}(k) \|^2  }_{\rm Term ~1}
  \\&  +\underbrace{   \left( g(k) -{1\over n} \sum_{i=1}^n \nabla f_i(x_i(k)) \right)^T(\bar{x}(k+1)-\bar{x}(k))}_{\rm Term ~2}
  \\&+\underbrace{ \bar{y}(k)^T ( \bar{x}(k+1)-\bar{x}(k) )   }_{\rm Term ~3}
  \\&  +\underbrace{ (\bar{x}(k)-\bar{x}(k+1) )^T {1\over n} \sum_{i=1}^n w_i(k)  } _{\rm Term ~4}.
\end{split}
\end{equation}
We  estimate the four terms on the right hand side of \eqref{bd-Fx}.

By  using \eqref{bd-barx}, we bound Term 1 as follows
\begin{equation}\label{bd-term1}
\begin{split}
{\rm Term ~1} & \leq  {\bar{\alpha}^2  L\over 2} \| \bar{y}(k)\|^2+ { c_1^2 L\over 2 n^2}    \|\tilde{y} (k)\|^2
\\& + { \bar{\alpha} c_1 L \over  n}   \| \bar{y}(k)\| \|\tilde{y} (k)\| .
\end{split}
\end{equation}

Note by the Jensen's inequality  that  for  $e=(e_1^T,\cdots, e_n^T)^T$,
$\left( \tfrac{  \sum_{i=1}^n \| e_i\|}{ n} \right)^2 \leq  {1\over n}  \sum_{i=1}^n \| e_i\|^2={1\over n} \|e\|^2 .$ Thus,
\begin{align}\label{bd-me}
 & \sum_{i=1}^n \| e_i\|  \leq \sqrt{n} \| e\|.
\end{align}
With \eqref{def-gk}  and Assumption \ref{ass-convex}(ii), we have
\begin{align*}
&\left \| g(k) -{1\over n} \sum_{i=1}^n \nabla f_i(x_i(k))  \right \|
  \leq  {1\over n} \sum_{i=1}^n \|  f_i (\bar{x}(k )) - f_i (x_i(k )) \|
 \\&  \leq {1\over n}  L \sum_{i=1}^n \|   \bar{x}(k )  -  x_i(k )\| \overset{ \eqref{bd-me}}{\leq } {L\over \sqrt{n}} \|\tilde{x}(k)\|.
\end{align*}   This combined with \eqref{bd-barx} produces that
\begin{align}\label{bd-term2}
  {\rm Term ~2} &\leq   { \bar{\alpha} L\over \sqrt{n}} \|\tilde{x}(k)\| \| \bar{y}(k)\|+   {  c_1 L\over n\sqrt{n}} \|\tilde{x}(k)\|    \|\tilde{y} (k)\|.
\end{align}

By using \eqref{recur-barx0}, we  bound Term 3 by the following
\begin{align}\label{bd-term3}
   {\rm Term ~3} &= \bar{y}(k)^T \Big( -\bar{\alpha}\bar{y}(k)- {(\mathbf{1}_n^T \bm{\alpha} - \bar{\alpha} \mathbf{1}_n^T)\otimes I_d \over n}  \tilde{y} (k)\Big) \notag
\\& \leq - \bar{\alpha} \| \bar{y}(k)\|^2+ {c_1 \over  n}   \| \bar{y}(k)\| \|\tilde{y} (k)\| .
\end{align}
By using \eqref{bd-barx}, \eqref{bd-me}, and $ab\leq {a^2\over 2c}+{cb^2\over 2},c>0$,  we  bound Term 4 as follows for any  $\mu \in (0,0.5],$
\begin{equation}\label{bd-term4}
\begin{split}
 {\rm Term ~4}&\leq \| \bar{x}(k)-\bar{x}(k+1)  \|   \sum_{i=1}^n  \|w_i(k) \|/n
 \\&\leq { \bar{\alpha} \over \sqrt{n}} \| \bar{y}(k)\|  \| \omega(k) \|+ {c_1 \over  n\sqrt{n}}   \| \bar{y}(k)\|   \| \omega(k) \|
 \\& \leq { \mu \bar{\alpha}^2 L \over 2} \| \bar{y}(k)\|^2+   {  1\over 2 \mu nL} \| \omega(k) \|^2
 \\& + { \mu\bar{\alpha}^2 L \over 2}    \| \bar{y}(k)\|^2+ {c_1^2 \over 2 n^3\mu \bar{\alpha}^2 L}   \| \omega(k) \|^2
 \\& =  \mu \bar{\alpha}^2 L   \| \bar{y}(k)\|^2 +a_1 \| \omega(k) \|^2,
\end{split}
\end{equation}
where  $ a_1\triangleq {   1 \over 2 \mu nL} + {c_1^2 \over 2 n^3\mu \bar{\alpha}^2 L}  .$

Therefore, by substituting \eqref{bd-term1}, \eqref{bd-term2}, \eqref{bd-term3}, and \eqref{bd-term4}    into  \eqref{bd-Fx}, we have that
 \begin{align}\label{bd-Fx2}
 &\mathbb{E}[v(k+1)|\mathcal{F}(k)] \leq v(k) -\left( \bar{\alpha} - (0.5+\mu) \bar{\alpha}^2 L  \right) \| \bar{y}(k)\|^2 \notag
 \\& +a_1 \| \omega(k) \|^2 + a_2 \|\tilde{x}(k)\| \| \bar{y}(k)\|+   a_3 \|\tilde{x}(k)\|    \|\tilde{y} (k)\|\notag
  \\& +   a_4   \|\tilde{y} (k)\|^2   +  a_5  \| \bar{y}(k)\| \|\tilde{y} (k)\|,
\end{align}
 where
  \begin{equation}\label{def-parameters}
\begin{split}
 &  a_2\triangleq { \bar{\alpha} L\over \sqrt{n}} ,  a_3\triangleq {  c_1 L\over n\sqrt{n}} , a_4\triangleq  { c_1^2 L\over  2 n^2} , a_5\triangleq {c_1 \over  n}(1+   \bar{\alpha} L).
\end{split}
\end{equation}
Taking the unconditional expectation on both sides of \eqref{bd-Fx2}  and summing it over $k$ from $0$ to $K,$ we get
 \begin{align}\label{bd-Fx3}
 & \mathbb{E}[v(k+1) ]\leq \mathbb{E}[v(0)]+a_1\sum_{k=0}^K  \mathbb{E}[\| \omega(k) \|^2 ]  \notag
 \\& -\left( \bar{\alpha} - (0.5+\mu) \bar{\alpha}^2 L  \right)\sum_{k=0}^K \mathbb{E}[ \| \bar{y}(k)\|^2]\notag
  \\&  + a_2\sum_{k=0}^K \mathbb{E}[\|\tilde{x}(k)\| \| \bar{y}(k)\|]+  a_3 \sum_{k=0}^K \mathbb{E}[\|\tilde{x}(k)\|    \|\tilde{y} (k)\|]
  \\& +  a_4 \sum_{k=0}^K \mathbb{E}[ \|\tilde{y} (k)\|^2]   +   a_5 \sum_{k=0}^K \mathbb{E}[\| \bar{y}(k)\| \|\tilde{y} (k)\| ]. \notag
\end{align}

Define $\bar{Y}_K \triangleq \left(\sum_{k=0}^K \mathbb{E}[\| \bar{y}(k)\|^2] \right)^{1/2}, \tilde{X}_K\triangleq \left(\sum_{k=0}^K \mathbb{E}[\|\tilde{x}(k)\|^2]\right)^{1/2},$  and $\tilde{Y}_K\triangleq \left(\sum_{k=0}^K \mathbb{E}[\|\tilde{y}(k)\|^2]\right)^{1/2}.$
Then by applying the  Cauchy-Schwarz inequality  to \eqref{bd-Fx3}, and using \eqref{bd-noisevr}, we obtain that
 \begin{align}\label{bd-Fx4}
\mathbb{E}[v(k+1) ] & \leq \mathbb{E}[v(0)]+ a_1 n\nu^2 \sum_{k=0}^K {1\over N(k)} \notag
 \\& - \left( \bar{\alpha} - (0.5+\mu) \bar{\alpha}^2 L  \right) \bar{Y}_K^2  + a_2\tilde{X}_K \bar{Y}_K \notag
  \\&+  a_3\tilde{X}_K\tilde{Y}_K +  a_4\tilde{Y}_K^2   +   a_5 \tilde{Y}_K\bar{Y}_K.
\end{align}

Recall from Lemma \ref{lem-exp} and $(a+b)^2\leq 2(a^2+b^2)$ that
\begin{align}
\tilde{Y}_K & \leq q_0 +  q_1 \Phi_K + q_2\bar{Y}_K, \label{bd-YK}\\
\tilde{X}_K & \leq q_0+   q_1 \Phi_K  + q_2\bar{Y}_K,\label{bd-XK}
\end{align}
where $q_0\triangleq \left(3   \over 1-\rho_2^2\right)^{1/2} \left(\mathbb{E}  [\| \tilde{x} (0 )  \|^2+ \| \tilde{y} (0) \|^2] \right)^{1/2},$
 $q_1\triangleq {  \sqrt{3  }  \over 1-\rho_2},$ $ q_2\triangleq  { \sqrt{3( c_1^2+   c_2^2L^2)}\over  1-\rho_2} ,$
and \begin{align}\label{def-phi}
 \Phi_K\triangleq \nu \Big(  n   \sum_{s=0}^{K }  \big(  N(k+1) ^{-1/2}   + N(k) ^{-1/2}\big)^2\Big)^{1/2}. \end{align}
This implies that
\begin{align*}
& \tilde{X}_K \bar{Y}_K  \leq  q_2\bar{Y}_K^2+ (q_0+   q_1 \Phi_K)\bar{Y}_K,
\\&\tilde{X}_K\tilde{Y}_K\leq   q_2^2\bar{Y}_K^2+2q_2(q_0+   q_1 \Phi_K)\bar{Y}_K+(q_0+   q_1 \Phi_K)^2
\\&\tilde{Y}_K^2 \leq  q_2^2\bar{Y}_K^2+2q_2(q_0+   q_1 \Phi_K)\bar{Y}_K+(q_0+   q_1 \Phi_K)^2,
\\&\tilde{Y}_K\bar{Y}_K\leq   q_2\bar{Y}_K^2+ (q_0+   q_1 \Phi_K)\bar{Y}_K   .
\end{align*}

 This incorporated with \eqref{bd-Fx4} produces
 \begin{align}\label{bd-Fx5}
 \mathbb{E}[v(k+1) ]& \leq \mathbb{E}[v(0)]+ a_1 n\nu^2\sum_{k=0}^K {1\over N(k)} \notag
  \\&  - a_0\bar{Y}_K^2 + b_0 \bar{Y}_K + c_0.
\end{align}
where $b_0\triangleq (a_2+a_5)(q_0+   q_1 \Phi_K)+2q_2(a_3+a_4)(q_0+   q_1 \Phi_K)$, $c_0\triangleq (a_3+a_4)(q_0+   q_1 \Phi_K)^2,$
and
\begin{align}\label{def-a0} a_0\triangleq \bar{\alpha} - (0.5+\mu) \bar{\alpha}^2 L -(a_2q_2+a_3q_2^2+a_4q_2^2+a_5q_2) .
 \end{align}

 Note from the definitions of $c_1$ and $c_2$ in Lemma \ref{lem-cond} that
\begin{align*}
q_2 &={\sqrt{3( c_1^2+   c_2^2L^2)}\over  1-\rho_2} %={3( c_1^2(1+L^2)+n\bar{\alpha}^2 )\over  1-\rho_2}
   ={\sqrt{3(1+L^2) \sum_{i=1}^n \alpha_i^2 -3n\bar{\alpha}^2 }\over  1-\rho_2}
\\& < {\sqrt{3n(1+L^2) } \alpha_{\max}  \over  1-\rho_2}  .
\end{align*}
This combined with \eqref{def-parameters} and $\mu \in (0,0.5]$ implies that
\begin{align*}
a_0 &>\bar{\alpha} -   \bar{\alpha}^2 L -  \bar{\alpha} \alpha_{\max} L \left(1+  d_{\alpha}  (  1/L+   \bar{\alpha}  ) \right)  {\sqrt{3 (1+L^2) } \over  1-\rho_2} \notag
\\& - \bar{\alpha}L  \Big(  d_{\alpha}     + {  d_{\alpha}^2 \bar{\alpha} \over  2   }\Big) { 3 (1+L^2)   \alpha_{\max}^2    \over  (1-\rho_2)^2}   {\rm~with~} d_{\alpha}= {c_1\over \sqrt{n} \bar{\alpha}}.
\end{align*}
This combined with  $-\bar{\alpha} \geq -\alpha_{\max}$ implies that
\begin{align*}
a_0&>\bar{\alpha} \Bigg(1
 - \alpha_{\max}^2 L   \Big(  d_{\alpha}     + {  d_{\alpha}^2  \alpha_{\max} \over  2   }\Big) { 3 (1+L^2)  \over  (1-\rho_2)^2} \\& -    \alpha_{\max}  L \Big(1+ \left(1+  d_{\alpha}  (  L^{-1}+   \alpha_{\max}  ) \right)  {\sqrt{3 (1+L^2) } \over  1-\rho_2}\Big)\Bigg) .
\end{align*}
It can be   seen that $a_0>0 $  for sufficiently small $\alpha_i>0.$

Since  $v(k)\geq 0,$    from \eqref{bd-Fx5} it follows that   for any $K\geq 1,$ \begin{align*}
 \mathbb{E}[v(0)]+ a_1 n\nu^2\sum_{k=0}^K {1\over N(k)}  - a_0\bar{Y}_K^2 + b_0 \bar{Y}_K + c_0 \geq 0.
\end{align*}
By recalling  the condition $\sum_{k=0}^{\infty} {1\over N(k)}<\infty$, we conclude form
\eqref{def-phi} that $\Phi_K<\infty$  for any $K\geq 1.$ Hence $b_0<\infty, c_0<\infty$,
and  $\bar{Y}_K$ is  uniformly bounded.
Since   $\{\bar{Y}_K\}$ is  an increasing sequence, we  conclude  that
$$\lim_{K\to \infty}\bar{Y}_K=\bar{Y}_{\infty}<\infty.$$
Similarly, from \eqref{bd-YK} and \eqref{bd-XK},  we  obtain  that
$$\lim_{K\to \infty} \tilde{Y}_K=\tilde{Y}_{\infty}<\infty {\rm~and~} \lim_{K\to \infty} \tilde{X}_K=\tilde{X}_{\infty}<\infty.$$

The above two equations imply  that
\begin{equation}\label{as-xy}
\begin{split}
& \sum_{k=0}^{\infty} \mathbb{E}[\| \bar{y}(k)\|^2] <\infty \Rightarrow \sum_{k=0}^{\infty} \| \bar{y}(k)\|^2 <\infty,\quad a.s.,
 \\& \sum_{k=0}^{\infty} \mathbb{E}[\|  \tilde{x}(k)\|^2] <\infty \Rightarrow \sum_{k=0}^{\infty} \|  \tilde{x}(k)\|^2 <\infty,\quad a.s.,
\\& \sum_{k=0}^{\infty} \mathbb{E}[\|  \tilde{y}(k)\|^2] <\infty\Rightarrow \sum_{k=0}^{\infty} \| \tilde{y}(k)\|^2 <\infty,\quad a.s.~.
 \end{split}
 \end{equation}
 Recall from \eqref{bd-noisevr} and  $\sum_{k=0}^{\infty} {1\over N(k)}<\infty$   that
\begin{align}\label{as-noise}
\sum_{k=0}^{\infty}  \mathbb{E}[\|w(k)\|^2]<\infty \Rightarrow \sum_{k=0}^{\infty} \| w(k)\|^2 <\infty,\quad a.s.~.
\end{align}
By using \eqref{bd-Fx2} and $ab\leq {a^2+b^2\over 2}$,  we obtain that
 \begin{align}\label{bd-cond}
 &\mathbb{E}[v(k+1)|\mathcal{F}(k)] \leq v(k) -\bar{\alpha}\left( 1- (0.5+\mu) \bar{\alpha}  L  \right) \| \bar{y}(k)\|^2 \notag
 \\& +a_1 \| \omega(k) \|^2 + { a_2 +a_3\over 2}   \|\tilde{x}(k)\|^2 +   \notag
  \\& +   {a_3+2a_4+a_5\over 2}   \|\tilde{y} (k)\|^2   +  { a_2+a_5 \over 2}   \| \bar{y}(k)\|^2.
\end{align}

We then use \eqref{as-xy}, \eqref{as-noise},  apply    the Robbins-Siegmund theorem in \cite{robbins1971convergence} to \eqref{bd-cond},   and conclude that $v(k)=F(\bar{x}(k )) -F^*$ converges  almost surely to some finite random variable.
Therefore, $\{\bar{x}(k )\}$ is almost surely  bounded.  Let $\bar{x} $ be a cluster point of $\{\bar{x}(k )\}$.
  Then there   exists a subsequence $k_s$  such that
$   \lim\limits_{s \rightarrow \infty  } \bar{x}(k_s )= \bar{x} .$
Note from \eqref{as-xy} that $ \lim\limits_{k \rightarrow \infty  } \| \tilde{x}(k)\|=0,~a.s.$,  and hence
\begin{align}  \label{limit-xik}
  \lim\limits_{s \rightarrow \infty  } x_i(k_s )= \bar{x},\quad \forall i\in \mathcal{V} .
\end{align}
Note by \eqref{as-xy} and \eqref{as-noise}   that $ \lim\limits_{k \rightarrow \infty  }\bar{y}(k)=0,~a.s. $ and $ \lim\limits_{k \rightarrow \infty  } w(k)=0,~a.s. ,$ respectively.   Then by  \eqref{recur-bary1} there holds
$ \lim\limits_{k \rightarrow \infty  } {1\over n} \sum_{i=1}^n f_i (x_i(k ))=0 . $
This incorporated with \eqref{limit-xik} produces $ {1\over n} \sum_{i=1}^n f_i (\bar{x})=0 $. Hence,
$\bar{x}$ is an optimal solution to the problem \eqref{problem1}.  Therefore,
$$\lim_{k \rightarrow \infty  }F(\bar{x}(k )) -F^* =F(\bar{x} ) -F^*=0,\quad a.s.$$
This completes the proof.
\hfill $\Box$

\subsection{Proof of Corollary \ref{cor1}}\label{app:A4}

Denote by $\beta\triangleq \alpha L.$  By  \eqref{def-alpha1} and $c_0<1$, we have that   $$\beta<{ c_0+1+2\sqrt{3}L-\sqrt{(c_0+1+2\sqrt{3}L)^2-4c_0} \over 2 }<c_0. $$
This implies that $\beta^2-(c_0+1+ 2\sqrt{3  }  L)\beta +c_0>0,$ and  hence
$(1  -    \beta)\left( c_0- \beta     \right) >    2\sqrt{3  }\beta L .$
By the definition of $\rho_2$ in Lemma \ref{lem-exp}, using  $\beta<c_0$ and
$c_0={ (1-\rho_1)^2\over (2-\rho_1) }$,
there holds $$1-\rho_2 >1-\rho_1-{\beta \over 2} -{\sqrt{c_0^2+4c_0} \over 2}
={c_0-\beta \over 2}.$$ Therefore,
\begin{align*}  &    (1  -    \beta)(  1-\rho_2 ) -    \sqrt{3  }\beta L
 >   (1  -    \beta){c_0-\beta \over 2}-    \sqrt{3  }\beta L >0.
\end{align*}
This implies that
 \begin{align}\label{bd-beta}  1  -    \beta  >{ \sqrt{3  }\beta L \over    1-\rho_2  } .\end{align}

Since all the agents take the same stepsize    $\alpha_i\equiv  \alpha $, we have
 $\bar{\alpha}=\alpha,~c_1=0,~c_2=\sqrt{n} \alpha,$  and
$ q_2= {\sqrt{3n }  \beta  \over  1-\rho_2} .$
Then by \eqref{def-parameters}, we obtain that $  a_2={  \beta\over \sqrt{n}} $  and $  a_3= a_4= a_5=0.$
Thus,  from \eqref{def-a0} it follows that
$ a_0= \alpha \left(1 - (0.5+\mu)  \beta -\tfrac{  \sqrt{3} \beta L } {  1-\rho_2}\right).$
This combined with \eqref{bd-beta} and $\mu\in (0,0.5)$ makes  $ a_0>0$.
The rest of the proof is the same as that of   Theorem  \ref{thm-as}.
\hfill $\Box$
\section{Proofs of Section  \ref{sec:sc}} \label{app:B}

\subsection{  Proof of Lemma \ref{lem1}.}\label{app:B1}
{\em   We first give an upper bound on $\| \bar{x}(k+1) -x^*\| $.}
By using     \eqref{recur-barx0},    \eqref{recur-bary1}, and $F(x )= {1\over n } \sum_{i=1}^n f_i(x)$, we obtain   that
\begin{align}\label{recur-barx2}
&\| \bar{x}(k+1) -x^*\|  \notag
\\& = \| \bar{x}(k)  -\bar{\alpha}\bar{y}(k)- {(\mathbf{1}_n^T \bm{\alpha} - \bar{\alpha} \mathbf{1}_n^T)\otimes I_d \over n}  \tilde{y} (k)-x^*\| \notag
\\&=\Big \| \bar{x}(k) -x^* - \bar{\alpha}  \nabla F(\bar{x}(k))  + {\bar{\alpha} \over n} \sum_{i=1}^n  \nabla f_i(\bar{x}(k))  \notag
\\&
\quad- {\bar{\alpha} \over n} \sum_{i=1}^n  \nabla f_i(x_i(k))- {\bar{\alpha} \over n} \sum_{i=1}^n w_i(k) \notag\\
 &\quad-{(\mathbf{1}_n^T \bm{\alpha} - \bar{\alpha} \mathbf{1}_n^T)\otimes I_d \over n}  \tilde{y} (k)\Big \|
 \notag \\& \overset{(a)}{ \leq  }  \left  \| \bar{x}(k) -\bar{\alpha}  \nabla F(\bar{x}(k)) -x^* \right\|
 \notag\\&\quad +{ \bar{\alpha} \over n} \left \|  \sum_{i=1}^n  \nabla f_i(\bar{x}(k)) -  \sum_{i=1}^n  \nabla f_i(x_i(k))\right\|
\\& \quad +  \bar{\alpha} \left\|  {  1\over n} \sum_{i=1}^n  \nabla f_i(x_i(k)) - \bar{y}(k)\right \|
+{c_1 \over n} \| \tilde{y} (k)\|
\notag \\&  \overset{(b)}{ \leq  }   \left  \| \bar{x}(k) -\bar{\alpha}  \nabla F(\bar{x}(k)) -x^* \right\|
\notag \\&  + { \bar{\alpha } L  \over n} \sum_{i=1}^n \| x_i(k) -  \bar{x}(k)  \|  +    {\bar{\alpha }\over n}\sum_{i=1}^n \left\| w_i (k)\right \|+{ c_1\over n} \| \tilde{y} (k)\|,\notag
\end{align}
where in (a)   we use the triangle
inequality and $\|\mathbf{1}_n^T \bm{\alpha} - \bar{\alpha} \mathbf{1}_n^T\|=c_1$,  and in (b) we use
 Assumption \ref{ass-convex}(ii).

We  introduce an inequality from \cite[Eqn. (2.1.24)]{nesterov2013introductory} on   the
$\eta$-strongly convex and $L$-smooth   function  $f(x)$,
\begin{align}\label{inequ-f}
 & (x-y)^T(\nabla f(x) -\nabla f(y)) \geq {\eta L \over \eta+L} \|x-y\|^2   \notag \\&
+{1\over \eta+L}  \|\nabla f(x)- \nabla f(y)\|^2,\quad \forall x,y\in \mathbb{R}^d.
\end{align}
By  $\alpha_i  \in (0, {2\over \eta+L}]$, we have that  $\bar{\alpha}\in (0, {2\over \eta+L}] $ and  ${2\over\bar{\alpha}}-\eta\geq L.$ Define $L' \triangleq {2\over \bar{\alpha}}-\eta$.
 \red{From Assumptions \ref{ass-convex}(ii) and \ref{ass-fun} it is seen that   the    function $F(x)={1\over n} \sum_{i=1}^n f_i(x)$ is  $\eta$-strongly convex and $L'$-smooth.}
Thus, by applying     \eqref{inequ-f} with $x=x(k)$ and $y=x^*,$ from  $ \nabla F(x^*)=0$  and ${2\over \eta+L'}=\bar{\alpha} $   it follows that
\begin{align*}
 & \left \| \bar{x}(k) -\bar{\alpha}  \nabla F(\bar{x}(k)) -x^* \right\|^2   \\&   =    \left \| \bar{x}(k)  -x^* \right\|^2+ \bar{\alpha}^2 \left \|   \nabla F(\bar{x}(k))  \right\|^2
\\&\quad  -2\bar{\alpha}  (\bar{x}(k)   -x^*)^T (\nabla F(\bar{x}(k)) -\nabla F(x^*))
\\& \leq  \left \| \bar{x}(k)  -x^* \right\|^2 +\bar{\alpha}^2 \left \|   \nabla F(\bar{x}(k))  \right\|^2
\\&\quad - 2\bar{ \alpha} \left({\eta L' \over \eta+L'} \|x_k-x^*\|^2+{1\over \eta+L'}  \|\nabla F(x_k) \|^2\right)
\\&  \leq  \Big(1-{2 \bar{\alpha}\eta L' \over \eta+L'} \Big) \left \| \bar{x}(k)  -x^* \right\|^2
  - \bar{\alpha} \Big({2\over \eta+L'}-\bar{\alpha} \Big)  \|\nabla F(x_k) \|^2
\\& \leq   \left(1-{2 \bar{\alpha}\eta L' \over \eta+L'} \right) \left \| \bar{x}(k)  -x^* \right\|^2
=\left(1- \bar{\alpha}\eta \right)^2  \left \| \bar{x}(k)  -x^* \right\|^2.
\end{align*} Then we can bound the first term  of     \eqref{recur-barx2}  by
 \begin{align}\label{bd-xstar}
& \left \| \bar{x}(k) -\alpha  \nabla F(\bar{x}(k)) -x^* \right\|  \leq (1-\bar{\alpha }\eta) \left \| \bar{x}(k)  -x^* \right\|.
\end{align}
  Therefore, by  plugging    \eqref{bd-xstar}  into    \eqref{recur-barx2} and using  \eqref{bd-me}, we have that
\begin{align}
&\| \bar{x}(k+1) -x^*\|   \leq  (1-\bar{\alpha }\eta)\| \bar{x}(k)-x^*\|  \notag\\&
 \qquad + {\bar{\alpha}L \over \sqrt{n}}  \| \tilde{x}(k)   \| +    {\bar{\alpha }\over \sqrt{n}} \left\| w (k)\right \|+{c_1 \over n} \| \tilde{y} (k)\|. \label{recur-barx21}
\end{align}

{\em Next, we   provide an  upper bound on $\|    \bar{y}(k)  \| .$}
 By    using $\sum_{i=1}^n \nabla f_i(x^*)=0  $    and \eqref{recur-bary1}, we obtain that
\begin{align*}
&\|    \bar{y}(k)  \|=\left \|  {1\over n} \sum_{i=1}^n (  w_i(k) + \nabla f_i(x_i(k))-\nabla f_i(x^*) )\right \|  \notag
\\&  \overset {(a)}{ \leq} {1\over n} \sum_{i=1}^n  \left \|  w_i(k) \right \|  + { L\over n} \sum_{i=1}^n  \| x_i(k) - x^*\|
\\&   \overset{\eqref{bd-me}} {\leq} {1\over \sqrt{n}} \left \|   w(k) \right \|  + { L\over \sqrt{n}}   \| x(k) - (\mathbf{1}_n \otimes I_d)   x^*\|
\notag \\ &
\overset{(b)}{  \leq}{1\over \sqrt{n}} \left \|   w(k) \right \|
 + { L\over \sqrt{n}}   \| x(k) - (\mathbf{1}_n \otimes I_d)   \bar{x}(k)\|
\\&\quad+{ L\over \sqrt{n}}   \| (\mathbf{1}_n \otimes I_d)   \bar{x}(k) - (\mathbf{1}_n \otimes I_d)   x^*\|  \\ &
   \overset{\eqref{def-tidlex}}{= }{1\over \sqrt{n}} \left \|   w(k) \right \|  + { L\over \sqrt{n}}   \| \tilde{x}(k)  \| +  L   \|   \bar{x}(k) -   x^*\|,
\end{align*}
where in (a)   we used the triangle  inequality and  Assumption \ref{ass-convex}(ii), in (b)  we added and  subtracted  the term $(\mathbf{1} \otimes I_d)   \bar{x}(k) $ and applied the triangle inequality.
  This combined with  \eqref{recur-tildex}  and \eqref{recur-yerr2}     produces
\begin{align}
 &\mathbb{E}[\| \tilde{x} (k+1)  \| | \mathcal{F}(k)] \leq  (\rho_1+ c_1L/ \sqrt{n} ) \| \tilde{x}(k)\|\notag
  \\& +{c_1\over \sqrt{n}} \left \|   w(k) \right \|  +    c_1L   \|   \bar{x}(k) -   x^*\|+ \alpha_{\max} \| \tilde{y}(k)\|,\label{recur-tildex2}
\end{align}
and
\begin{align}  \label{recur-yerr3}
 & \mathbb{E}[\| \tilde{y} (k+1)  \| | \mathcal{F}(k)] \leq     ( \rho_1 +\alpha_{\max}L) \| \tilde{y}(k)\| \notag
\\& +   (L+ c_2L^2/\sqrt{n}) \| \tilde{x}(k) \|      + \mathbb{E}[ \| w(k+1)-w(k ) \| | \mathcal{F}(k)]\notag
\\& +{ c_2 L \over \sqrt{n}} \left \|   w(k) \right \| + c_2 L^2   \|   \bar{x}(k) -   x^*\|.
\end{align}

Therefore,   by taking unconditional expectations on both sides of    \eqref{recur-barx21}, \eqref{recur-tildex2},   and \eqref{recur-yerr3},   we prove the lemma.
\hfill $\Box$

\subsection{Proof of Theorem \ref{thm1}} \label{app:B2}
 The spectral radius of the nonnegative matrix $J(\bm{\alpha})$ satisfying  $\rho(J(\bm{\alpha}))<1$ is equivalent to
 that all leading principle minors of $I_3-J(\bm{\alpha})$ are positive.
 Let $J_1=1-\bar{\alpha }\eta $. Then $det(I_1-J_1)=\bar{\alpha }\eta>0.$
Let $$J_2= \begin{pmatrix}
1- \bar{\alpha }\eta  &  {\bar{ \alpha}L \over \sqrt{n}}
 \\  c_1 L&   \rho_1+ {c_1L \over \sqrt{n} }   \end{pmatrix} .$$
  Then by $d_{\alpha}={c_1 \over \sqrt{n} \bar{\alpha}} $ and   $\kappa= L/\eta$, there holds
   \begin{equation}\label{bd-detJ2}
   \begin{split}det(I_2-J_2)&= det \begin{pmatrix}
 \bar{\alpha }\eta &- {\bar{ \alpha}L \over \sqrt{n}}
 \\ -c_1 L&  1-(\rho_1+ {c_1L \over \sqrt{n} })   \end{pmatrix}
%  \\&=\bar{\alpha } \eta( 1- \rho_1) -{c_1L\bar{\alpha }\eta \over \sqrt{n} } -{\bar{ \alpha} c_1L^2 \over \sqrt{n}}
% \\&=\bar{\alpha }\eta( 1- \rho_1) - {     d_{\alpha} \bar{\alpha }^2L(L+\eta)  }
 \\&= \bar{\alpha }L \big(  1- \rho_1  -     \bar{\alpha } d_{\alpha} \kappa(L+\eta) \big)/\kappa .
 \end{split}
 \end{equation}

Note that \begin{align*} &det(I_3- J(\bm{\alpha}))
\\& = det
 \begin{pmatrix}
 \bar{\alpha }\eta   & -{\bar{ \alpha}L \over \sqrt{n}}  & -{c_1 \over n}
 \\ -c_1 L&  1-(\rho_1+ {c_1L \over \sqrt{n} }) & -\alpha_{\max}  \\
- c_2 L^2 &-(L+{ c_2L^2 \over \sqrt{n}}) &1-(\rho_1 +\alpha_{\max}L )\end{pmatrix}
 \\&=   - c_2 L^2 det
 \begin{pmatrix}
  -{\bar{ \alpha}L \over \sqrt{n}}  & -{c_1 \over n}
 \\  1-(\rho_1+ {c_1L \over \sqrt{n} }) & -\alpha_{\max}   \end{pmatrix}
 \\&\quad +(L+{ c_2L^2 \over \sqrt{n}}) det
 \begin{pmatrix}
 \bar{\alpha }\eta \  &  -{c_1 \over n}
 \\ -c_1 L  & -\alpha_{\max}   \end{pmatrix}
 \\&\quad+(1-(\rho_1 +\alpha_{\max}L ))det
 \begin{pmatrix}
 \bar{\alpha }\eta   & -{\bar{ \alpha}L \over \sqrt{n}}
 \\ -c_1 L&  1-(\rho_1+ {c_1L \over \sqrt{n} })     \end{pmatrix}
 \\&= - c_2 L^2\left( {\bar{ \alpha}L \alpha_{\max} \over \sqrt{n}} + {c_1 \over n} \left(1- \rho_1- {c_1L \over \sqrt{n} }\right)   \right)
 \\& \quad-\left(L+{ c_2L^2 \over \sqrt{n}}\right)  \left(\alpha_{\max}  \bar{\alpha }\eta \   +{c_1^2 L \over n}\right)
 \\&\quad+(1-(\rho_1 +\alpha_{\max}L ))  \bar{\alpha }L \big(  1- \rho_1  -     \bar{\alpha } d_{\alpha} \kappa(L+\eta) \big)/\kappa.
 \end{align*}
Note by $d_{\alpha}= {\sqrt{\sum_{i=1}^n (\alpha_i-\bar{\alpha})^2} \over \sqrt{n} \bar{\alpha}} $
 and \eqref{def-alpha}  that
\begin{align}\label{def-c1c2}
 c_1=\sqrt{n} \bar{\alpha} d_{\alpha},~ c_2=\sqrt{c_1^2+n\bar{\alpha} ^2}=\sqrt{n} \bar{\alpha} \sqrt{d_{\alpha}^2+1}.
\end{align}
 Therefore,  \begin{align}   \label{bd-detJ3}
  &det(I_3- J(\bm{\alpha}))\notag
 %\\&= - \sqrt{n} \bar{\alpha} \sqrt{d_{\alpha}^2+1} L^2\left( {\bar{ \alpha}L \alpha_{\max} \over \sqrt{n}} + { \bar{\alpha} d_{\alpha} \over \sqrt{n}} \left(1- \rho_1-  \bar{\alpha}d_{\alpha}L  \right)   \right)\notag
% \\& -\left(L+   \bar{\alpha} \sqrt{d_{\alpha}^2+1}L^2  \right)  (\alpha_{\max}  \bar{\alpha }L/\kappa   +  \bar{\alpha}^2 d_{\alpha}^2  L )\notag
% \\&+ \bar{\alpha }L(1-(\rho_1 +\alpha_{\max}L ))\left( (1- \rho_1)/\kappa - \bar{\alpha}d_{\alpha} (L+\eta)  \right)\notag
 \\&= -  \bar{\alpha}^2 L^2\sqrt{d_{\alpha}^2+1}\left(  \alpha_{\max} L +     d_{\alpha}   \left(1- \rho_1-  \bar{\alpha} d_{\alpha} L \right)   \right)\notag
 \\& -\bar{\alpha } L^2  \left(1+   \bar{\alpha}L \sqrt{d_{\alpha}^2+1}  \right)  (\alpha_{\max}  + \kappa\bar{\alpha}d_{\alpha}^2    )/ \kappa \notag
 \\&+ \bar{\alpha }L(1- \rho_1 -\alpha_{\max}L ) \left(  1- \rho_1  -   (\kappa+1)\bar{\alpha}L d_{\alpha} \right) /\kappa \notag
 \\&=\bar{\alpha }L \Big(  (1- \rho_1 -\alpha_{\max}L ) \left(  1- \rho_1  -   (\kappa+1)\bar{\alpha}L d_{\alpha} \right) /\kappa \notag
 \\& \qquad -  L  \left(1+   \bar{\alpha}L \sqrt{d_{\alpha}^2+1}  \right)  (\alpha_{\max}  + \kappa\bar{\alpha}d_{\alpha}^2    )/ \kappa \notag
 \\&\qquad -  \bar{\alpha} L \sqrt{d_{\alpha}^2+1}\left(  \alpha_{\max} L +     d_{\alpha}   \left(1- \rho_1-  \bar{\alpha} d_{\alpha} L \right)   \right)\Big). \notag
 \\&\geq \bar{\alpha }L \Big(  (1- \rho_1)^2/\kappa  -(\alpha_{\max}L  + \bar{\alpha}L(\kappa+1) d_{\alpha}  )(  1- \rho_1 ) /\kappa \notag
 \\& \qquad -  L  \left(1+   \bar{\alpha}L \sqrt{d_{\alpha}^2+1}  \right)  (\alpha_{\max}  + \kappa\bar{\alpha}d_{\alpha}^2    )/ \kappa
 \\&\qquad -  \bar{\alpha} L \sqrt{d_{\alpha}^2+1}\left(  \alpha_{\max} L +     d_{\alpha}   (1- \rho_1)   \right)\Big). \notag
 \end{align}

  From \eqref{bd-stepsize} it is seen that $ \alpha_iL<\beta^*$, hence  $ \alpha_{\max}L<\beta^*$ and
 $ \bar{\alpha}L<\beta^*$.   Then  by  \eqref{bd-detJ3}  and \eqref{bd-stepsize}, we have that
 \begin{equation*}
 \begin{split}
   & det(I_3- J(\bm{\alpha}))
   \\ & \geq \bar{\alpha }L \Big(  (1- \rho_1)^2  -\beta^*(1 +  (\kappa+1) d_{\alpha}  )(  1- \rho_1 )
 \\& \qquad -  \beta^* \left(1+  \beta^* \sqrt{d_{\alpha}^2+1}  \right)  (1 + \kappa d_{\alpha}^2    )
 \\&\qquad -  \beta^*\kappa \sqrt{d_{\alpha}^2+1}\left( \beta^*  +     d_{\alpha}   (1- \rho_1)   \right) \Big)/\kappa
 \\&= \bar{\alpha }L \big(  -c_3(\beta^* )^2-c_4\beta^* +(1- \rho_1)^2  \big)/\kappa>0.
 \end{split}
 \end{equation*}
 Note from \eqref{bd-stepsize} that  $\bar{\alpha}L <   \tfrac{1- \rho_1}{ d_{\alpha} \kappa(L+\eta) }.$
 Hence by \eqref{bd-detJ2} it is   seen that  $det(I_2-J_2)>0$.
 Therefore, \eqref{bd-stepsize} is a sufficient condition for  guaranteeing $\rho(J(\bm{\alpha})) <1.$

From $N(k)=\lceil q^{-2k} \rceil$  and  \eqref{bd-noise} it follows that
$ \mathbb{E}[\| w (k) \|] \leq  \sqrt{n} \nu q^k $ for any $ k\geq 0.$
 Then by    using \eqref{recursion-z0} and the triangle equality, we obtain that
\begin{align}\label{recur-z}
&z(k+1)  \leq  J(\bm{\alpha})z(k)  +\begin{pmatrix}
 \bar{ \alpha} \nu  \\ c_1\nu \\  \sqrt{n} \nu (1+q) +c_2 L   \nu \end{pmatrix}  q^k\notag  \\
&\leq J(\bm{\alpha})^{k+1}z(0)  + \sum_{t=0}^k  J(\bm{\alpha})^t q^{k-t} \begin{pmatrix}
 \bar{ \alpha} \nu  \\ c_1 \nu \\  \sqrt{n} \nu (1+q) +c_2   L   \nu
\end{pmatrix} .
\end{align}
\red{By noting that   $ \rho(J(\bm{\alpha}))<1$, $J(\bm{\alpha})^k$ converges to zero at the linear rate $\mathcal{O}( \rho(J(\bm{\alpha}))^k)$ (see \cite[Eqn. (7.10.5)]{Meyer00}).}  Thus,
\begin{equation}\label{bd-errz}
\begin{split}
&z(k)   \leq \mathcal{O}(\rho(J(\bm{\alpha}))^{k}) +  \sum_{t=0}^{k-1} \mathcal{O}(\rho(J(\bm{\alpha}))^t)    q^{k-1-t}.
\end{split}
\end{equation}

  Note that  for any $\rho<q$,
\begin{align*}
\sum_{t=0}^{k}  \rho ^t    q^{k-t}  =q^{k}  \sum_{t=0}^{k}  (\rho/q) ^t  \leq {q^{k} \over 1-\rho/q}= {q^{k+1} \over q-\rho } ,
\end{align*}
while for  any $\rho>q$,
$\sum_{t=0}^k  \rho ^t    q^{k-t}     \leq   \tfrac{q^{k+1} }{  \rho-q } .$
\red{Hence  $\sum_{t=0}^{k-1}  \rho ^t    q^{k-t}     \leq   { \max\{p,q\}^k  \over  | \rho-q |}  .$
This combined with \eqref{bd-errz}  proves  the theorem.}
\hfill $\Box$

\subsection{Proof of Corollary \ref{cor2}}\label{app:B3}

From   $\alpha_i\equiv \alpha$,   \eqref{def-z}, and \eqref{def-hatJ}  it follows  that $d_{\alpha}=0 $,  $\alpha_{\max}=\alpha $, and  $J(\bm{\alpha})=\hat{J}( \alpha )$.
Thus,   $det(I_2-J_2)>0$ by  \eqref{bd-detJ2}.
Then by defining $\beta\triangleq \alpha L$ and using   \eqref{bd-detJ3}, we obtain by  \eqref{stepsize}  that
\begin{align} \label{det}
&det(I_3- \hat{J}( \alpha ))\notag
  \\& =-{\beta \over \kappa } \left((1+\kappa) \beta^2 +  (2-\rho_1)\beta -(1- \rho_1 )^2\right) >0. \end{align}
Therefore, the step-size  \eqref{stepsize} makes  $\rho(\hat{J}( \alpha )) <1.$

\red{Note by  $\rho(\hat{J}( \alpha ))<q$ that $\rho( \hat{J}( \alpha )/q)<1$. Hence
\begin{equation*}
\begin{split}
&\sum_{t=0}^k \hat{J}( \alpha )^t q^{k-t}= q^k\sum_{t=0}^k  \big(\hat{J}( \alpha )/q \big)^t
\\&= q^k \big(I_3-  \hat{J}( \alpha )/q \big)^{-1}\big(I_3- (\hat{J}( \alpha )/q)^{k+1}\big) .
\end{split}
\end{equation*}
This  together with \eqref{recur-z} and $J(\bm{\alpha})=\hat{J}( \alpha )$  implies that for sufficiently  large  $k$,
\begin{equation*}
\begin{split}
&z(k+1)  \approx  \hat{J}( \alpha ) ^{k+1} z(0) \\& + q^k \big(I_3-  \hat{J}( \alpha ) /q \big)^{-1} \begin{pmatrix}
 \bar{ \alpha} \nu  \\ c_1 \nu \\  \sqrt{n} \nu (1+q) +c_2   L   \nu
\end{pmatrix}   .
\end{split}
\end{equation*}
Then the result follows by  $\bar{\alpha}=\alpha, c_1=0, $ and $c_2=\sqrt{n}\alpha $.}
\hfill $\Box$

\blue{\subsection{Proof of Corollary \ref{cor3}}\label{app:B4}
%
%$\begin{pmatrix}
%  \alpha  \eta  & { - \alpha L \over \sqrt{n}}  & 0
% \\ 0&  1-\rho_1  & -\alpha \\
% -\sqrt{n}\alpha L^2 &-L-\alpha L^2  &1-\rho_1 -\alpha L \end{pmatrix} $
%
From $N(k)=B$  and  \eqref{bd-noise} it follows that
$ \mathbb{E}[\| w (k) \|] \leq  \sqrt{n/B} \nu .$
 Similarly to \eqref{recur-z}, we obtain that
 \begin{align*}
&z(k+1)  \leq  J(\bm{\alpha}) ^{k+1}z(0)  + \sum_{t=0}^k  { J(\bm{\alpha})^t \over \sqrt{B}} \begin{pmatrix}
 \bar{ \alpha} \nu  \\ c_1 \nu \\  2\sqrt{n} \nu  +c_2   L   \nu
\end{pmatrix} .
\end{align*}
It has been shown in Corollary \ref{cor2} that  the step-size satisfying  \eqref{stepsize} makes
 $\bar{\alpha}=\alpha, c_1=0, $ and $c_2=\sqrt{n}\alpha $, $J(\bm{\alpha})=\hat{J}( \alpha )$, and
$\rho(\hat{J}( \alpha )) <1.$ Hence
  \begin{align*}
&z(k+1)  \leq  \hat{J}( \alpha )^{k+1}z(0)  + \nu\sum_{t=0}^k  { \hat{J}( \alpha )^t \over \sqrt{B}} \begin{pmatrix}
\alpha  \\ 0 \nu \\  2\sqrt{n}   +\sqrt{n}\alpha  L
\end{pmatrix} .
\end{align*}
Note by \cite[Eqn. (7.10.11)]{Meyer00}) that
$ \sum_{p=0}^{\infty}    \hat{J}(\alpha)^p=(I_3-\hat{J}(\alpha))^{-1}$.
 Therefore, $z_k$ converge to  $\limsup_{k\to \infty} z_k$
 with    a geometric  rate $\mathcal{O}\big( \rho(\hat{J}(\alpha)) ^k\big)$, and
 \begin{align*}
&\limsup_{k\to \infty}z(k)  \leq { \nu (I_3-\hat{J}(\alpha))^{-1}\over \sqrt{B}}  \begin{pmatrix}
   \alpha    \\  0 \\  (2 + \alpha   L  ) \sqrt{n}
\end{pmatrix} .
\end{align*}
By $(I_3-\hat{J}(\alpha))^{-1}={(I_3-\hat{J}(\alpha))^{*} \over det( I_3-\hat{J}(\alpha)} $,  we have that
\begin{align*}
&\limsup_{k\to \infty} \mathbb{E}[\| \bar{x}(k) -x^*\| ]\leq  { \nu \over \sqrt{B} det( I_3-\hat{J}(\alpha)} \times
\\& \big[\alpha(
 (1-\rho_1)( 1-\rho_1-\alpha L)-\alpha L(1+  \alpha L  ))+  (2 + \alpha   L  )     \alpha^2 L \big]
 \\&= {\alpha \nu ( (1-\rho_1)^2+\rho_1\alpha L) \over \sqrt{B} det( I_3-\hat{J}(\alpha))}  ,\quad {\rm and ~}
\\& \limsup_{k\to \infty} \mathbb{E}[ \| x(k) - (\mathbf{1}_n \otimes I_d)\bar{x}(k)  \| ]
\\& \leq {\nu   \over \sqrt{B} det( I_3-\hat{J}(\alpha)} \big[
\alpha \sqrt{n}\alpha^2 L^2  + \alpha^2\eta (2 + \alpha   L  )  \sqrt{n}    \big]
\\& = {\alpha^2 \sqrt{n} \nu  (\alpha L^2  +  \eta (2 + \alpha   L  )     ) \over \sqrt{B} det( I_3-\hat{J}(\alpha)}.
 \end{align*}
  This combined with \eqref{det} proves the results,
 \hfill $\Box$

\subsection{Proof of Theorem \ref{thm4}}\label{app:B5}
Since $N(k)=\lceil (k+1)^{2\theta} \rceil  $ and  $ \rho(J(\bm{\alpha}))<1$,  by using the similar procedures for deriving
\eqref{recur-z} and \eqref{bd-errz},  we have
\begin{align*}
z(k )&  \leq  J(\bm{\alpha})z(k)  +\begin{pmatrix}
 \bar{ \alpha} \nu  \\ c_1\nu \\  2\sqrt{n} \nu   +c_2 L   \nu \end{pmatrix}  k^{-\theta}
\\&\leq \mathcal{O}(\rho(J(\bm{\alpha}))^{k}) +  \sum_{t=1}^{k } \mathcal{O}(\rho(J(\bm{\alpha}))^{k-t})   t^{-\theta}  =\mathcal{O}( k^{-\theta} ).
\end{align*}
Thus, $z(k)\leq C_3 k^{-\theta} $ for some $C_3>0$, and hence $\| z(k )\| \leq   \epsilon  $ for any $  k\geq  K_3(\epsilon) =  \left({C_3 \over \epsilon}\right)^ {1 /\theta }$.
Thus, the  number of   communication rounds required is $ 2|\mathcal{E}|K_3(\epsilon) $, and the  number of sampled gradients required is bounded by
\begin{align*}
&\sum_{k=0}^{K_3(\epsilon)} N(k)=\sum_{k=0}^{K_3(\epsilon)} (k+1)^{2\theta}  \leq
 \int_{1}^{ K_3(\epsilon)+1  } t^{2\theta} dt
\\&= \tfrac{t^{2\theta+1} }{ 2\theta}\Big |_{1}^{ K_3(\epsilon) }
 =  (2\theta+1)^{-1} \left(\tfrac{ C_3}{\epsilon }\right)^{(2\theta+1)/\theta}  . \end{align*}
This completes the proof. \hfill $\Box$}
%\bibliographystyle{IEEEtran}
%\bibliography{DSO}

% Generated by IEEEtran.bst, version: 1.13 (2008/09/30)
\def\cprime{$'$}

 \end{document}